\documentclass[sn-mathphys-num]{sn-jnl}

\usepackage{graphicx}%
\usepackage{multirow}%
\usepackage{amsmath,amssymb,amsfonts}%
\usepackage{amsthm}%
\usepackage{mathrsfs}%
\usepackage[title]{appendix}%
\usepackage{xcolor}%
\usepackage{textcomp}%
\usepackage{manyfoot}%
\usepackage{booktabs}%
\usepackage{algorithm}%
\usepackage{algorithmicx}%
\usepackage{algpseudocode}%
\usepackage{listings}%

\usepackage{float}
\usepackage{amsbsy, mathrsfs}
\usepackage{bbm}
\usepackage[normalem]{ulem}
\usepackage{caption}
\usepackage{subcaption}

\usepackage{tikz}
\usetikzlibrary{positioning}
\usetikzlibrary{calc}
\usetikzlibrary{arrows}
\usetikzlibrary{decorations.pathmorphing,decorations.markings}
\usetikzlibrary{shapes}
\usetikzlibrary{patterns}

\definecolor{blue}{rgb}{0.23,0.58,0.89}

\tikzset{
  pics/dist_arc/.style args={#1,#2,#3,#4}{
     code={
       \draw[thick, #4] (-0.25,0.25) -- (0.25,0.25);
       \draw[thick, #4] (-0.5,-0.04) -- (-0.15,0.25);
       \draw[thick, #4] (0.5,-0.04) -- (0.15,0.25);
       \draw[->, thick, #4] (0,0.5) -- (0,0.25);
       \node[#3] (#1) at (0,0) {#2};
     }
  },
  minimum size=2.5em
}

\begin{document}

\title{Solving decision problems with endogenous uncertainty and conditional information revelation using influence diagrams}


\author[1]{\fnm{Olli} \sur{Herrala}}\email{olli.herrala@aalto.fi}
\author[2]{\fnm{Tommi} \sur{Ekholm}}\email{tommi.ekholm@fmi.fi}
\author*[1]{\fnm{Fabricio} \sur{Oliveira}}\email{fabricio.oliveira@aalto.fi}

\affil*[1]{\orgdiv{Department of Mathematics and Systems Analysis}, \orgname{Aalto University}, \orgaddress{\street{PO Box 11000}, \city{Espoo}, \postcode{FI-00076 AALTO}, \country{Finland}}}

\affil[2]{\orgname{Finnish Meteorological Institute},\orgaddress{\city{Helsinki}, \country{Finland}}}


\abstract{Mathematical programming formulations of influence diagrams can bridge the gap between representing and solving decision problems. However, they suffer from both modeling and computational limitations. Aiming to address modeling limitations, we show how to incorporate conditionally observed information within the mathematical programming representation of the influence diagram. Multi-stage stochastic programming models use conditional non-anticipativity constraints to represent such uncertainties, and we show how such constraints can be incorporated into the influence diagram formulations. This allows us to consider the two main types of endogenous uncertainty simultaneously, namely decision-dependent information structure and decision-dependent probability distribution. 

Additionally, we apply a subdiagram decomposition to improve both computational efficiency and modeling capabilities. Under suitable conditions, this decomposition allows for considering continuous decision variables arising from, e.g., investment sizing decisions, leading to better solutions than a discretization of the continuous decisions. Finally, our proposed framework is illustrated with a large-scale cost-benefit problem regarding climate change mitigation, simultaneously considering technological research and development, and optimal emission trajectories. }




\maketitle

\section{Introduction}

A common assumption in decision making under uncertainty is that the stochastic processes, particularly the state probabilities and/or observed values, are not influenced by the previously made decisions. Thus, the uncertainty is said to be exogenous. Here, we focus on a much less explored class of stochastic problems presenting endogenous uncertainty. In this more general setting, decisions can affect the uncertainty faced in later stages through \emph{decision-dependent probabilities} or \emph{conditionally observed information} \citep{hellemo2018}, which makes these problems challenging from both a modeling and computational standpoint. To solve such problems, we employ a combination of influence diagrams and mixed-integer programming.

Influence diagrams \citep{howard2005influence} are one of the common ways to represent decision problems. Their strength lies in conveying the problem in an intuitive, visual form, representing the various interdependencies between decisions, chance events, and outcomes simultaneously. A notable benefit of using influence diagrams is that modeling decision-dependent probabilities is straightforward. Additionally, \citet{Salo2022} argue that conditionally observed information can also be represented through influence diagrams. They present a simple problem to demonstrate this, but their discussion is brief and, for example, does not analyze the implications of such endogenous uncertainty on computational performance.

Despite their capability to represent endogenous uncertainty, influence diagrams are otherwise limited in their ability to represent decision problems by two restricting requirements, namely \textit{regularity} (decisions are made sequentially in a fixed order) and \textit{no-forgetting} (all previously made decisions and observed uncertainties are remembered when making a decision) \citep{Shachter1986}. To circumvent these requirements, thus broadening our modeling capabilities, this paper considers limited-memory influence diagrams (LIMID) from \citet{lauritzen2001representing}, that is, influence diagrams that do not satisfy regularity and no-forgetting requirements. Their main example is a partially observed Markov decision process (POMDP) without memory, that is, the decisions in each period are made based on an imperfect observation with no memory of past decisions or observations. Limited-memory influence diagrams have been used in a variety of contexts, such as airplane maintenance planning \citep{cai2009influence}, cancer treatment strategies \citep{van2007selecting} and power plant monitoring \citep{agogino2013integrating}, highlighting their flexibility in modeling decision processes. For the remainder of this paper, \textit{influence diagrams} are assumed to contain LIMIDs, following the convention of \citet{koller2009probabilistic}.

The solution methods for influence diagrams have relied on reducing the diagram by adding and removing arcs and merging nodes \citep{howard2005influence, maua2012solving}, possibly combined with a subdiagram decomposition \citep{lee2021submodel}. However, these approaches suffer from both modeling and computational issues \citep{bielza2000structural, bielza2011review}. One of the modeling issues for influence diagrams is that these solution approaches assume an expected utility maximization problem with no additional constraints. This precludes the use of, e.g., chance or budget constraints and (conditional) value-at-risk to model the decision maker's risk preferences. 

Seeking to overcome the limitations of the previous solution methods, two mixed-integer programming (MIP) formulations for (limited-memory) influence diagrams have been recently introduced: rooted junction trees \citep{parmentier2020integer} and decision programming \citep{Salo2022}, which leverage the advancements in solving stochastic programming (SP) models. Compared to the existing solution methods for influence diagrams, these formulations have advantages including the possibility to add different risk constraints, being able to solve the models using an off-the-shelf solver, and applying, e.g., decomposition methods to make larger models computationally tractable. Of these two formulations, the rooted junction tree model of \citet{parmentier2020integer} generally results in better computational performance \citep{herrala2023RJT}. 

Although virtually non-existent within the influence diagram setting, tackling conditionally observed information has been discussed extensively in the stochastic programming context \citep[e.g.][]{apap2017}. However, decision-dependent probabilities pose significant modeling challenges in this setting, requiring specific dependency structures and problem-specific reformulations \citep{hellemo2018}. 

To fill the gap between these two domains, our paper explores the possibility of incorporating conditionally observed information in mixed-integer formulations for IDs in two ways: 1) by using observation nodes in the influence diagram; and 2) by implementing conditional non-anticipativity constraints \citep{apap2017}, originally used in multistage stochastic programming, within the rooted junction tree model \citep{parmentier2020integer}. This results in a framework capable of addressing both types of endogenous uncertainty simultaneously, while still retaining prior computationally favorable properties of the mathematical model such as linearity, convexity and the ability to incorporate, e.g., chance constraints \citep{herrala2023RJT}.

We additionally explore the possibility of introducing continuous variables to the problem. The formulation presented in \citet{parmentier2020integer} accommodates only discrete decisions, and \citet{Salo2022} acknowledge the limitations of using influence diagram-based models for problems involving continuous decisions, a challenge discussed in more detail in \citet{bielza2011review}. Decisions such as order or production quantities and investment amounts arise in various problem settings, and we show that these decision spaces can be approximated with discrete sets of alternatives, at the expense of trading off solution quality and computational performance. To circumvent this issue, we show that if the problem has a separable structure, the submodel-tree decomposition proposed by \citet{lee2021submodel} can be employed, making it possible to incorporate continuous decisions, as long as they do not affect the endogenous probabilities in the model. In our computational experiments, the solution times using this decomposition are similar to those using a very coarse discretization of the decision space. 

This paper is structured as follows. In Section \ref{section:uncertainty}, we present an overview of stochastic programming and the influence diagram formulation by \citet{parmentier2020integer}. In Section \ref{section:contributions}, our methodological contributions are described in detail, starting from conditionally observed information and continuing with the role of subdiagram decomposition in incorporating continuous variables. The computational performance of these contributions is then assessed in Section \ref{section:experiments}. In Section \ref{section:score} we illustrate the use of the framework by considering a larger-scale problem of climate change cost-benefit analysis. Section \ref{section:conclusions} concludes and provides directions for further development.

\section{Modeling problems with endogenous and exogenous uncertainties}
\label{section:uncertainty}

\subsection{Decision-making under endogenous uncertainty}

Before discussing the use of influence diagrams for modeling decision problems and how to obtain optimal strategies using mixed-integer formulations based on rooted junction trees, we first give an overview of the application and research areas from which our research builds upon. In mathematical programming, stochastic programming (SP) is one of the most widespread frameworks for decision-making under uncertainty. In general, SP casts decision problems subject to parametric uncertainty as deterministic equivalents in the form of large-scale linear or mixed-integer linear programming (LP/MILP) models that can be solved with standard optimization techniques. SP models often assume that the uncertainty is \emph{exogenous}, i.e., independent of the decisions made in previous stages. This is methodologically convenient, for the deterministic equivalent model has the same nature as its stochastic counterpart, retaining important characteristics such as linearity, or more generally, convexity.

Solution approaches for multi-stage stochastic programming (MSSP) under exogenous uncertainty are often based on formulating the deterministic equivalent problem using a scenario tree, as described in, e.g., \citet{Ruszczynski1997}. A scenario tree represents the structure of the uncertain decision process, and non-anticipativity constraints (NACs) \citep{rockafellar1991scenarios} are employed to enforce the information structure in the formulation. NACs state that a decision must be the same for two scenarios if those scenarios are indistinguishable when making the decision. 

In the more general setting of \emph{endogenous} uncertainty, the decisions made at previous stages can affect the uncertainty faced in later stages. \cite{hellemo2018} propose a taxonomy that classifies endogenous uncertainties into two distinct types. In Type 1 problems, earlier decisions influence the later events' probability distribution (i.e., realizations and/or the probabilities associated with each realization). For example, deciding to perform maintenance on a car engine influences the probability of it breaking down in the future. In Type 2 problems, the \emph{information structure} is influenced by the decision-making. Continuing with the car example, deciding to inspect the engine does not affect the probability of breakdown, but provides information that enables a better-informed maintenance decision. Type 3 combines the Type 1 and Type 2 endogenous uncertainties. 

Type 2 endogenous uncertainty has been more widely addressed in the SP literature, perhaps due to one of its subclasses having a strong connection to exogenously uncertain problems. In the taxonomy presented in \citet{hellemo2018}, this specific type of endogenous uncertainty is called conditional information revelation. For the sake of terminology consistency, we refer to this as \emph{conditionally observed information}. In this subclass, the decisions only affect the time at which the (exogenous) uncertainty is revealed to the decision maker. One of the earliest publications on such uncertainty is \cite{jonsbraaten1998}, where the authors describe a branching algorithm for solving a subcontracting problem. \citet{Goel2006} consider a process network problem where the yield of a new process is uncertain prior to installation. Other applications include open pit mining \citep{boland2008multistage}, clinical trial planning for drug development \citep{colvin2010modeling} and technology project portfolio management \citep{solak2008stochastic}. 

Similar solution methods are employed in problems with exogenous uncertainty and conditionally observed information. The main difference is that conditional observation requires the use of \emph{conditional non-anticipativity constraints} (C-NACs), as the distinguishability between scenarios at any given stage is dependent on earlier decisions. This conditional dependency results in disjunctive constraints that require specific reformulation techniques \citep{apap2017}. The main challenge arising from this approach is that the number of constraints rapidly increases with problem size, resulting in computational intractability for large problems. \citet{apap2017} also propose omitting redundant constraints, in an attempt to mitigate the tractability issues. In their example considering a production planning problem, this results in roughly a 99\% decrease in the problem size. Despite these substantial improvements, the reduced model is still considerably large and cannot be solved to the optimum within a reasonable computation time under their experimental setting, illustrating how challenging such problems are.

Type 1 endogenous uncertainty is more challenging from a mathematical modeling standpoint because the uncertain events themselves depend on earlier decisions. Consequently, a scenario tree-based representation cannot be posed, as the scenario probabilities in a scenario tree cannot depend on decisions. Therefore, the well-established solution techniques for MSSP cannot be directly applied to these problems. Despite these challenges, some discussion on Type 1 endogenous uncertainty is found in the literature. \citet{Peeta2010} discuss the fortification of a structure in a network, where the probability of failure depends on the fortification decision. \citet{dupacova2006optimization} presents a summary of problems with Type 1 uncertainties and \citet{escudero2020some} present solution approaches to multi-stage problems where the first-stage decisions influence the scenario probabilities in later stages. Examples of such problems can be found in \citet{zhou2022stochastic} and \citet{li2023strategic}. Finally, reformulations and custom algorithms for various Type 1 problems are further summarized in \citet{hellemo2018}. However, these approaches assume specific relationships between decisions and probabilities and result in non-convex nonlinear formulations. Consequently, these approaches are not easily generalizable to different problems. 

\subsection{Influence diagrams}
\label{subsec:id}

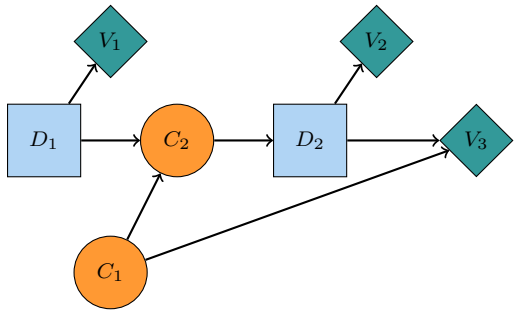
\begin{figure}[H]
\centering
\begin{tikzpicture}
    [decision/.style={fill=blue!40, draw, minimum size=2.5em, inner sep=2pt}, 
    chance/.style={circle, fill=orange!80, draw, minimum size=2.5em, inner sep=2pt},
    value/.style={diamond, fill=teal!80, draw, minimum size=2.5em, inner sep=2pt},
    optimization/.style={ellipse, fill=teal!80, draw, minimum size=2em, inner sep=2pt},
    scale=1.75, font=\scriptsize]
    
    \node[chance]   (L) at (-0.25, 1.5)  {$I_{pre}$};
    \node[chance]   (L1) at (0.75, 2.25)  {$R_1$};
    \node[chance]   (L2) at (0.75, 0.75)  {$R_2$};
    \node[decision] (1) at (1.75, 2.25)  {$A_1$};
    \node[decision] (2) at (1.75, 0.75)  {$A_2$};
    \node[chance]   (F) at (2.75, 1.5)  {$I_{post}$};
    \node[value]    (T) at (3.75, 1.5)  {$U_{int}$};   
    \node[value]    (V1) at (2.75, 2.5)  {$V_1$};   
    \node[value]    (V2) at (2.75, 0.5)  {$V_2$};     
    \draw[->, thick] (L) -- (L1);
    \draw[->, thick] (L) -- (L2);
    \draw[->, thick] (L) -- (F);
    \draw[->, thick] (L1) -- (1);
    \draw[->, thick] (L2) -- (2);
    \draw[->, thick] (1) -- (F);
    \draw[->, thick] (2) -- (F);
    \draw[->, thick] (F) -- (T);
    \draw[->, thick] (1) -- (V1);
    \draw[->, thick] (2) -- (V2);
    \draw[->, thick] (L) to[out=60,in=135,distance=1.5cm] (V1);
    \draw[->, thick] (L) to[out=-60,in=-135,distance=1.5cm] (V2);
\end{tikzpicture}
\caption{An influence diagram representation of a decision problem}
\label{fig:nmon}
\end{figure}

As discussed earlier, modeling decision-dependent probabilities in stochastic programming poses considerable challenges. To circumvent this limitation, we rely instead on \emph{influence diagrams}, which allow for an intuitive and structured representation of how uncertainties and decisions are mutually dependent. An influence diagram $G_{ID} = (N,A_{ID})$ is an acyclic graph formed by nodes $j \in N = N^C \cup N^D \cup N^V$ and arcs $a \in A_{ID} = \{(i,j) \mid i,j \in N\}$\footnote{We use the subscript $ID$ for $G_{ID}$ and $A_{ID}$ to distinguish between influence diagrams and rooted junction trees. Rooted junction trees are described in Section \ref{subsec:rjt}.}. Nodes $N^C$ and $N^D$ are the sets of chance and decision nodes, respectively, and $N^V$ is a collection of value nodes representing the consequences incurred from the decisions made at nodes $N^D$ and the chance events realized at nodes $N^C$. In Fig. \ref{fig:nmon}, the decision nodes are represented by squares, the chance nodes by circles, and the value nodes by diamonds. 

To illustrate the use of influence diagrams, we use the N-monitoring problem in Fig. \ref{fig:nmon}, introduced in \citet{Salo2022}. To ease the upcoming discussion in this paper, let us assume the following problem setting: a company is launching a new product, and the product has an initial public interest level represented by the chance node $I_{pre}$. There are two sales representatives in different areas, each obtaining a report $R_i$ about the interest level by surveying their locations. Note that the number of areas could be extended to an arbitrary number N, hence the name for the problem. The report is considered to be an uncertain estimate of reality, represented as a chance node. Based on their report, each representative decides (decision nodes $A_i$) whether or not to sell the product in their area, resulting in a profit or loss represented by the value node $V_i$. In this problem, neither regularity nor no-forgetting assumption hold, as the decisions are made in parallel with no communication between the decision makers, making the diagram in Fig. \ref{fig:nmon} a LIMID. 

Introducing a product in a given area results in more overall visibility for the product on, e.g., social media. This results in an influence between $A_i$ and the (uncertain) post-decision state of the overall public interest level $I_{post}$, for which the realization translates into a utility value represented by the value node $U_{int}$. This node represents the utility function of the company, i.e., how valuable the different outcomes in $I_{post}$ are to the company, allowing us to combine the monetary profit/loss in nodes $V_i$ with the monetary value the company sees in different outcomes of $I_{post}$. This represents a more general problem of moving from a pre-decision state $I_{pre}$ to a stochastic post-decision state $I_{post}$ by making parallel decisions informed by imperfect reports of the pre-decision state. In \cite{Salo2022}, the N-monitoring problem represents a safety-critical system where different sensors monitor the state of the system, and actions must be taken based on these sensor readings to keep the system operational. 

Each decision and chance node $j \in N^C \cup N^D$ can assume a state $s_j$ from a discrete and finite set of states $S_j$. For a decision node $j \in N^D$, $S_j$ represents the available choices; for a chance node $j \in N^C$, $S_j$ is the set of possible realizations. For instance, the possible levels of public interest represented in nodes $I_{pre}$ and $I_{post}$ form a discrete set, such as $S_j = \{$``low'', ``medium'', ``high''$\}$ for $j \in \{I_{pre}, I_{post}\}$. The arcs  $(i,j)$ in $A_{ID} = \{(i,j) \mid i,j \in N\}$ represent influence between nodes. In Fig. \ref{fig:nmon}, the arcs are represented by arrows between the nodes. Before defining this notion of influence further, let us first define a few necessary concepts.

The \emph{information set} comprises the immediate predecessors of a given node $j \in N$ and is defined as $I(j) = \{i \in N \mid (i,j) \in A_{ID}\}$. In the graphical representation, this corresponds to the set of nodes that have an arrow pointing directly to node $j$. For example, in Fig. \ref{fig:nmon}, the information set of $I_{post}$ consists of $I_{pre}$, $A_1$ and $A_2$, that is, the initial level of interest ($I_{pre}$) along with the actions taken by the regional representatives ($A_i$) affect the probabilities of different levels of post-decision public interest ($I_{post}$). The decision $s_j \in S_j$ made in each decision node $j \in N^D$ and the conditional probabilities of the states $s_j \in S_j$ in each chance node $j \in N^C$ depend on their \emph{information state} $s_{I(j)} \in S_{I(j)}$, where $S_{I(j)} = \prod_{i \in I(j)} S_i$. Referring to our example, the probabilities of different outcomes in $I_{post}$ are conditional on the decisions in $A_1$ and $A_2$ and the random outcome in $I_{pre}$. Let us define $X_j \in S_j$ as the realized state at a chance node $j \in N^C$. Using the notion of information states, the conditional probability of observing a given state $s_j$ for $j \in N^C$ is $\mathbb{P}(X_j = s_j \mid X_{I(j)} = s_{I(j)})$. 

 For a decision node $j \in N^D$, let $Z_j: S_{I(j)} \to S_j$ be a mapping between each information state $s_{I(j)} \in S_{I(j)}$ and decision $s_j \in S_j$. That is, $Z_j(s_{I(j)})$ defines a \emph{local decision strategy}, which represents the choice of some $s_j \in S_j $ in $j \in N^D$, given the information $s_{I(j)}$. In the N-monitoring example in Fig. \ref{fig:nmon}, the decisions $A_i$ are made knowing the outcome of the report in $R_i$. Note that we do not consider mixed strategies, where each information state would be mapped to an arbitrary probability distribution over $S_j$. Instead, we only consider deterministic strategies that can be represented by an indicator function $\mathbb{I}: S_{I(j)} \times S_j \to \{0,1\}$ defined so that
\begin{align}
    \mathbb{I}(s_{I(j)}, s_j) = \begin{cases} 1, &\text{ if } Z_j \text{ maps } s_{I(j)} \text{ to } s_j \text{, i.e., } Z_j(s_{I(j)}) = s_j; \\ 0, &\text{ otherwise.} \end{cases} \label{eq:strategy_indicator}
\end{align}
A (global) decision strategy is the collection of local decision strategies in all decision nodes: $Z = (Z_j)_{j \in N^D}$, selected from the set of all possible strategies $\mathbb{Z}$.

At each value node $v \in N^V$, a real-valued utility function $U_v : S_{I(v)} \to \mathbbm{R}$ maps the information state $s_{I(v)}$ of $v$ to a utility value $U_v$. In the N-monitoring example, the profit/loss made in area $i$ depends on the initial interest level in $I_{pre}$ and the decision $A_i$ of whether to sell the product in the area. The default objective is to maximize the expected utility of a strategy, but other objectives such as conditional Value-at-Risk can also be used \citep{Salo2022,herrala2023RJT}.

\subsection{Rooted junction trees}
\label{subsec:rjt}
As shown in \citet{Salo2022}, it is possible to obtain a mixed-integer linear programming (MILP) model directly from the influence diagram representation of the problem. However, the authors observe that the model size increases exponentially with the number of nodes, resulting in computational challenges with relatively small problems. To mitigate this exponential growth, \citet{parmentier2020integer} proposes first reformulating the influence diagram into a \emph{rooted junction tree} (RJT) $G_{RJT} = (V, A_{RJT})$, a directed graph consisting of clusters $C \in V$ of nodes $j \in N$, and arcs between these clusters, with the underlying undirected graph (obtained by replacing the directed edges $A_{RJT}$ with undirected edges) being a tree. The core idea of the RJT reformulation is to construct the joint probability distributions corresponding to each cluster, which results in a smaller MIP model compared to the path-based approach in \citet{Salo2022}. The computational performance of the RJT reformulation is compared to Decision Programming \citep{Salo2022} in \citet{herrala2023RJT}, where the authors show that the growth in model size as more nodes are introduced to the diagram is significantly slower in RJT models. In this subsection, we first present key properties of a rooted junction tree before discussing the algorithms for converting an influence diagram to an RJT and the resulting MILP reformulation.

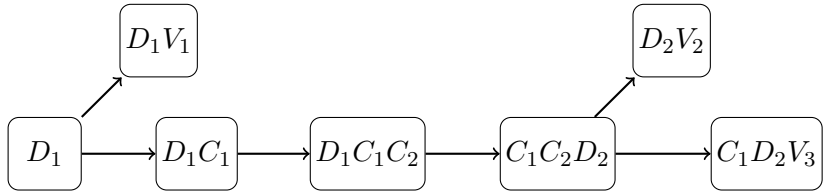
\begin{figure}[H]
\centering
\begin{tikzpicture}
    [cluster/.style={fill=white!80, draw, minimum size=2.5em, inner sep=2pt, rounded corners}]
     \node[cluster] (1) at (0, 0)       {$I_{pre}$};
     \node[cluster] (2) at (1.5, 0)     {$I_{pre} R_1$};
     \node[cluster] (3) at (3.25, 0)    {$I_{pre} R_1 A_1$};
     \node[cluster] (4) at (5.25, 0)    {$I_{pre} A_1 R_2$};
     \node[cluster] (5) at (7.5, 0)     {$I_{pre} A_1 R_2 A_2$};
     \node[cluster] (6) at (9.85, 0)    {$I_{pre} A_1 A_2 I_{post}$};
     \node[cluster] (7) at (12, 0)      {$I_{post} U_{int}$};
     \node[cluster] (8) at (4, -1.25)   {$I_{pre} A_1 V_1$};
     \node[cluster] (9) at (8.5, -1.25) {$I_{pre} A_2 V_2$};
     \draw[->, thick] (1) -- (2);
     \draw[->, thick] (2) -- (3);
     \draw[->, thick] (3) -- (4);
     \draw[->, thick] (4) -- (5);
     \draw[->, thick] (5) -- (6);
     \draw[->, thick] (6) -- (7);
     \draw[->, thick] (3) -- (8);
     \draw[->, thick] (5) -- (9);
\end{tikzpicture}
\caption{A gradual RJT representation of the ID in Fig. \ref{fig:nmon}}
\label{fig:nmon_rjt}
\end{figure}

The first important property of an RJT is the \emph{running intersection property}, i.e., if a node $j \in N$ is in two clusters of the tree, it is also in all clusters on the (undirected) path between these clusters. From this property, it follows that the subgraph of $G_{RJT}$ induced by a node $j$ (formed by clusters $C \in V$ for which $j \in C$, and the arcs connecting such clusters) is a rooted tree. For example, in Fig. \ref{fig:nmon_rjt}, the node $A_1$ appears in five clusters that form a rooted tree. 

More specifically, \citet{parmentier2020integer} consider gradual RJTs, where each cluster is the \emph{root cluster} of exactly one node $j \in N$, and the root cluster of node $j$ is denoted with $C_j$. Correspondingly, Fig. \ref{fig:nmon} and Fig \ref{fig:nmon_rjt} both have nine nodes or clusters, each cluster being the root node of the last node listed in the cluster. A root cluster is defined as the root of the subgraph of the RJT induced by a node $j$. For example, the root cluster of the subtree induced by $A_1$ is denoted as $C_{A_1}$. Finally, it is required that $I(j) \subset C_j$ for all $j \in N$. These properties result in a convenient structure where for each pair of adjacent clusters $(C_i,C_j) \in A_{RJT}$, we have $C_j \setminus C_i = j$ (as the only new node that can be introduced in $C_j$ without violating the running intersection property is $j$), and the joint probability distribution of all nodes in $C_j$ can thus always be obtained from the probability distribution of $C_i$ and the conditional probability $\mathbb{P}(X_j = s_j \mid X_{I(j)} = s_{I(j)})$ using the chain rule. As will be described next, this allows the problem to be formulated as a mixed-integer programming (MIP) model, which allows for employing standard techniques widely available in off-the-shelf solvers.

While any gradual RJT fulfilling these requirements can be used to formulate the MIP model \eqref{eq:rjt-obj}-\eqref{eq:rjt-vars}, large clusters result in larger models and the associated computational challenges. To avoid forming RJTs possessing large clusters, \citet{parmentier2020integer} present two algorithms for converting an influence diagram into an RJT with minimal clusters. In summary, both algorithms build the RJT by traversing the diagram according to a reverse topological order, and for each node $j$ in $N$, adds a cluster $C_j$, adds the node to any previously introduced clusters where the node is needed to satisfy the running intersection property and $I(j) \in C_j$, and finally adds the necessary arcs from $C_j$ to previously introduced clusters, repeating until the RJT has been built. The first algorithm uses a given topological ordering, while the second also forms a topological ordering resulting in smaller clusters.

To illustrate the conversion from an influence diagram to an RJT, we use the root cluster of node $R_2$ in Fig. \ref{fig:nmon_rjt} ($\{I_{pre}, A_1, R_2\}$) as an example. This cluster contains node $R_2$ along with its information set, allowing us to use the chain rule to determine the probability of the nodes in each cluster attaining specific states, as discussed later in conjunction with the MIP formulation \eqref{eq:rjt-obj}-\eqref{eq:rjt-vars}. In addition to nodes $R_2$ and $I(R_2)=\{I_{pre}\}$, the example cluster also contains $A_1$. While this node does not influence the probabilities in $R_2$, it is needed in the cluster $C_{I_{post}}$ to determine $\mathbb{P}(X_{I_{post}}=s_{I_{post}} \mid X_{I(I_{post})} = s_{I(I_{post})})$. Therefore, to maintain the running intersection property, $A_1$ must appear in all clusters on the path between $C_{A_1}$ and $C_{I_{post}}$. A more detailed description of the properties of gradual RJTs and the algorithms for building them is out of scope for this paper, and we instead refer the reader to \citet{parmentier2020integer} and \citet{herrala2023RJT}. 

Let us define the binary variable $z(s_j \mid s_{I(j)})$ that takes value 1 if $Z_j(s_{I(j)}) = s_j$, and 0 otherwise, for all $j \in N^D$, $s_j \in S_j$, and $s_{I(j)} \in S_{I(j)}$. These variables correspond to the indicator function \eqref{eq:strategy_indicator}, representing local decision strategies at each decision node $j \in N^D$. Additionally, we define variables $\mu_{C_j} \ge 0$ representing the joint probability distribution of nodes $i \in C_j$. The expected utility maximization problem corresponding to the gradual RJT can then be formulated as in \citet{parmentier2020integer}:
\begin{align}
    \max &\sum_{j \in N^V} \sum_{s_{C_j} \in S_{C_j}} \mu_{C_j}(s_{C_j}) u_{C_j}(s_{C_j}) \label{eq:rjt-obj}\\
    \text{s.t. } & \sum_{s_{C_j} \in S_{C_j}} \mu_{C_j}(s_{C_j}) = 1, \ \forall j \in N \label{eq:rjt-probsum}\\
    & \sum_{\substack{s_{C_i} \in S_{C_i}, \\ s_{C_i \cap C_j} = s^*_{C_i \cap C_j}}} \mu_{C_i}(s_{C_i})= \sum_{\substack{s_{C_j} \in S_{C_j}, \\ s_{C_i \cap C_j} = s^*_{C_i \cap C_j}}} \mu_{C_j}(s_{C_j}),\ \forall (C_i,C_j) \in A_{RJT}, s^*_{C_i \cap C_j} \in S_{C_i \cap C_j} \label{eq:rjt-localconsistency}\\
    & \mu_{C_j}(s_{C_j}) = \mu_{\overline{C}_j}(s_{\overline{C}_j}) \mathbb{P}(X_j=s_j \mid X_{I(j)}=s_{I(j)}), \ \forall j \in N^C \cup N^V, s_{C_j} \in S_{C_j} \label{eq:rjt-prob}\\
    & \mu_{C_j}(s_{C_j}) = \mu_{\overline{C}_j}(s_{\overline{C}_j})z(s_j \mid s_{I(j)}), \ \forall j \in N^D, s_{C_j} \in S_{C_j} \label{eq:rjt-dec}\\
    & \mu_{C_j}(s_{C_j}) \ge 0, \ \forall j \in N, s_{C_j} \in S_{C_j} \label{eq:rjt-mu}\\
    & z(s_j \mid s_{I(j)}) \in \{0,1\}, \ \forall j \in N^D, s_j \in S_j, s_{I(j)} \in S_{I(j)}. \label{eq:rjt-vars}
\end{align}

The objective function \eqref{eq:rjt-obj} is the expected utility associated with the strategy $Z \in \mathbb{Z}$ represented by the decision variables $z$. Following the definition of a gradual RJT, $I(j) \subset C_j$ and there is exactly one cluster corresponding to each value node $j \in N^V$. We thus set utility values associated with these clusters to $u_{C_j}(s_{C_j}) = U_j(s_{I(j)})$ using the utility function defined in Section \ref{subsec:id}. Constraints \eqref{eq:rjt-probsum} and \eqref{eq:rjt-mu} state that the decision variables $\mu_{C_j}$ must represent valid probability distributions, i.e., have nonzero probabilities that sum to one, and constraint \eqref{eq:rjt-localconsistency} enforces local consistency between adjacent clusters. Here, local consistency means that the distribution $\mu_{C_i \cap C_j}$ must be the same when obtained as a marginal distribution from $C_i$ or $C_j$. The states of the nodes in the intersection of clusters $C_i$ and $C_j$ are denoted with $s_{C_i \cap C_j}$, and a specific state combination in this set is denoted with $s^*_{C_i \cap C_j}$. 

Constraints \eqref{eq:rjt-prob} and \eqref{eq:rjt-dec} propagate the probability information in the junction tree using the chain rule $\mathbb{P}(X,Y)=\mathbb{P}(Y)\mathbb{P}(X|Y)$. For notational brevity, we use $\overline{C}_j = C_j \setminus j$, for which $\mu_{\overline{C}_j}(s_{\overline{C}_j}) = \sum_{s_j} \mu_{C_j}(s_{C_j})$. Combined with the local consistency constraints \eqref{eq:rjt-localconsistency}, this ensures that $\mu_{\overline{C}_j}$ can be used in the chain rule to obtain $\mu_{C_j}$. It should be noted that while constraint \eqref{eq:rjt-dec} contains a product of two decision variables, the decision strategy variables $z$ are binary, making \eqref{eq:rjt-dec} an indicator constraint. This allows one to linearize the product using methods discussed in, e.g., \citet{mitra1994tools}, and the problem can be considered an instance of mixed integer linear programming (MILP).

As an example, we discuss the variables and constraints for the first four clusters in Fig. \ref{fig:nmon_rjt}. For completeness, the corresponding parts of the model are presented in Appendix \ref{app:nmon}. For the first cluster, we introduce variables $\mu_{C_{I_{pre}}} \ge 0$ and add a constraint $\sum_{s_{I_{pre}} \in S_{I_{pre}}} \mu_{C_{I_{pre}}}(s_{I_{pre}}) = 1$. As there is only one node in the cluster, constraint \eqref{eq:rjt-prob} becomes simply $\mu_{C_{I_{pre}}}(s_{I_{pre}}) = \mathbb{P}(X_{I_{pre}}=s_{I_{pre}})$, that is, the probability distribution of the one-node root cluster must match the given distribution for that node. For the next cluster, we define $\mu_{C_{R_1}}$ in a similar way, but constraint \eqref{eq:rjt-prob} now becomes $\mu_{C_{R_1}}(s_{I_{pre}}, s_{R_1}) = \mu_{\overline{C}_{R_1}}(s_{I_{pre}})\mathbb{P}(X_{R_1}=s_{R_1} \mid X_{I_{pre}}=s_{I_{pre}})$, where $\mu_{\overline{C}_{R_1}}(s_{I_{pre}}) = \mu_{C_{I_{pre}}}(s_{I_{pre}})$ in this case, as $C_{R_1} \setminus R_1 = \{I_{pre}\} = C_{I_{pre}}$. Finally, the intersection $C_{I_{pre}} \cap C_{R_1}$ is $I_{pre}$ and constraint \eqref{eq:rjt-localconsistency} states that for each $s^*_{I_{pre}} \in S_{I_{pre}}$, we must have $\mu_{C_{I_{pre}}}(s^*_{I_{pre}}) = \sum_{s_{R_1}} \mu_{C_{R_1}}(s^*_{I_{pre}}, s_{R_1})$.

The process is otherwise the same for the second and third clusters, but for the third cluster, constraint \eqref{eq:rjt-dec} is used instead of \eqref{eq:rjt-prob}. For the fourth cluster $C_{V_1}$, we notice that $\overline{C}_{V_1}$ is not equal to the parent cluster $C_{A_1}$. Instead, constraint \eqref{eq:rjt-localconsistency} enforces local consistency and the distribution $\overline{C}_{V_1}$ is therefore effectively defined using the parent cluster distribution. Finally, as $V_1 \in N^V$ is a value node, this cluster contributes $\sum_{s_{I_{pre}} \in S_{I_{pre}}, s_{A_1} \in S_{A_1}, s_{V_1} \in S_{V_1}} \mu_{C_{V_1}}(s_{I_{pre}},s_{A_1},s_{V_1}) u_{V_1}(s_{I_{pre}},s_{A_1},s_{V_1})$ to the objective function.

As mentioned earlier, the main advantage of using a mixed-integer programming formulation of an influence diagram is the ability to incorporate risk aversion in the form of, e.g., chance constraints. Additionally, the assumptions on regularity and no-forgetting are violated by limited-memory influence diagrams such as the one in Fig. \ref{fig:nmon}, rendering traditional solution methods for influence diagrams unsuitable. Even for problems such as the finite memory POMDP example in \citet{lauritzen2001representing}, where well-established solution methods for similar problems (Markov decision processes) exist, the lack of memory and perfect observations makes the problem significantly harder to solve \citep{kara2023convergence}. In contrast, the employment of mixed-integer programming formulation to solve influence diagrams requires none of such assumptions.

\section{Methodological developments}
\label{section:contributions}
\subsection{Conditionally observed information}

As originally proposed, the approaches to formulate influence diagrams as MILP models most prominently focus on problems with Type 1 endogenous uncertainty, i.e., decision-dependent probabilities. However, many MSSP problems involve conditionally observed information and being able to model this within influence diagrams would make the models more generally applicable. After we present the concept of conditionally observed information, Section \ref{subsection:cnac_obs} will discuss the details of implementing this within the model \eqref{eq:rjt-obj}-\eqref{eq:rjt-vars}, and Section \ref{section:decomp} focuses on a decomposition approach for improving both the computational and modeling capabilities of the methods presented in this paper.

The previously discussed influence diagram in Fig. \ref{fig:nmon} represents the so-called N-monitoring problem \citep{Salo2022}, where a single underlying variable $I_{pre}$ is observed by $N$ independent decision makers (in Fig. \ref{fig:nmon}, $N=2$) through reports represented by $R_i$. This version of the model assumes that the reports are always available, but in reality, these reports might have a cost associated with their acquisition. Consequently, the cost may be higher than the value gained from the report and the optimal solution would thus be not to pay for the report. If we modify the problem so that the decision maker (DM) can choose whether or not to acquire the report, we say that observing the information in the report is conditional on the decision to pay for the report. We refer to this variant in Fig. \ref{fig:dist_arc_basic} as the conditional N-monitoring problem.

\begin{figure}[H]
\centering
\begin{tikzpicture}
    [decision/.style={fill=blue!40, draw, minimum size=2.5em, inner sep=2pt}, 
    chance/.style={circle, fill=orange!80, draw, minimum size=2.5em, inner sep=2pt},
    value/.style={diamond, fill=teal!80, draw, minimum size=2.5em, inner sep=2pt},
    optimization/.style={ellipse, fill=teal!80, draw, minimum size=2em, inner sep=2pt},
    scale=1.75, font=\scriptsize]
    
    \node[chance]   (L) at (-0.25, 1.5)  {$I_{pre}$};
    \node[chance]   (L1) at (0.75, 2.25)  {$R_1$};
    \node[decision] (1) at (2.25, 2.25)  {$A_1$};
    \node[chance]   (F) at (3.25, 1.5)  {$I_{post}$};
    \node[value]    (T) at (4.25, 1.5)  {$U_{int}$};   
    \node[value]    (V1) at (3.25, 2.5)  {$V_1$};      
    \node    (2) at (1.5, 0.75)  {...};    
    \node[decision]  (I) at (1.5, 3.5) {$D_1$};
    \draw[] (1.5,2.65) pic{dist_arc={int1, , ,}};
    \draw[thick] (L1) -- (int1.west);
    \draw[thick] (I) --  (int1.north);
    \draw[->, thick] (int1.east) --  (1);
    \draw[->, thick] (L) -- (L1);
    \draw[->, thick] (L) -- (F);
    \draw[->, thick] (1) -- (F);
    \draw[->, thick] (2) -- (F);
    \draw[->, thick] (F) -- (T);
    \draw[->, thick] (1) -- (V1);
    \draw[->, thick] (L) -- (2);
    \draw[->, thick] (L) to[out=60,in=135,distance=1.5cm] (V1);
\end{tikzpicture}
\caption{The conditional N-monitoring problem, where the report $R_1$ is conditionally observed when making the decision in $A_1$. The flow of information from $R_1$ to $A_1$ is conditional on $D_1$.}
\label{fig:dist_arc_basic}
\end{figure}

A key concept with conditionally observed information is \emph{distinguishability}. The outcome of a report itself, say $s_{R_i}$ or $s'_{R_i}$, is independent of the DM's choice to procure the report, but can only be observed if the DM pays for the report. Therefore, the DM cannot see the difference between $s_{R_i}$ and $s'_{R_i}$ when making the decision $A_i$ unless the report is procured. Otherwise, the report's outcome is not observed and the states $s_{R_i}$ and $s'_{R_i}$ are indistinguishable at $A_i$.

Our formulations for conditionally observed information focus on two key elements. First, we have the decisions or random events that the observation is conditional on, which we denote as the \emph{distinguishability set} $T_{i,j} \subset N^C \cup N^D$. Second, the observation depends on a \emph{distinguishability condition} $F_{i,j}: S_{T_{i,j}} \to \{0,1\}$. Here, $i \in N^C \cup N^D$ denotes the conditionally observed node and $j \in N^D$ is the decision node where that information is available if the distinguishability condition is fulfilled. 

In most cases, the conditionally observed node is a chance node, but it can also be a decision node. Such cases might arise in the context of distributed decision-making, where the decisions are made by multiple agents and observing the decision of another agent does not happen automatically, or problems in which previous decisions are not remembered by default and the decision maker must instead pay a price to retrieve information on past decisions. Distributed decision making in the context of influence diagrams is discussed in, e.g., \citet{detwarasiti2005influence}, and \citet{piccione1997interpretation} make a closely related point that decision makers can often affect what they remember by choosing to keep track of information (including decisions) they would otherwise forget.

Using the notion of distinguishability sets and conditions, we can define \emph{conditional arcs} \begin{equation}
    a_c \in A_c = \{(i,j) \mid i \in N^C \cup N^D, j \in N^D, \ T_{i,j} \neq \emptyset\}
\end{equation} 
to describe conditionally observed information in influence diagrams. Specifically, we say that a conditional arc $a_c$ from node $i \in N^C \cup N^D$ to node $j \in N^D$ is active (i.e., node $i$ is observed when making the decision corresponding to node $j$) if $F_{i,j}(s_{T_{i,j}}) = 1$. If $T_{i,j}$ is empty, there is no conditional observation of $i$ in $j$ and, thus, no conditional arc between these nodes exists. The concept is illustrated in Fig. \ref{fig:dist_arc_basic}. 
When making the decision in $A_i$, the report is available only if the decision maker chooses to acquire the report. Hence, the distinguishability set is $T_{R_i,A_i}=\{D_i\}$, and the distinguishability condition is $F_{R_i,A_i}(s_{D_i})=\mathbb{I}(s_{D_i} = \text{``acquire report"})$, where the indicator function $\mathbb{I}(\,\cdot\,)$ is defined as 
\begin{equation*}
    \mathbb{I}(x = x^*) = 
    \begin{cases} 
    1\,(true), &\text{if } x = x^*, \\
    0\,(false), &\text{otherwise.} 
    \end{cases}
\end{equation*}

If the distinguishability set includes more than one node, alternative functions $F_{i,j}$ might be employed for modeling the conditional dependencies between the nodes. For example, if there are several projects that reveal the same information in node $i \in N$ and completing any of these projects is sufficient for the information to be revealed, $F_{i,j}(s_{T_{i,j}}) = \bigvee_{k \in T_{i,j}} \mathbb{I}(s_k = s_k^*)$ can be used; or if all of the projects are required for the information revelation, $F_{i,j}(s_{T_{i,j}}) = \bigwedge_{k \in T_{i,j}} \mathbb{I}(s_k = s_k^*)$ is appropriate. An example of such conditions is found in \citet{tarhan2009stochastic}, where different uncertainties in oil field development are gradually revealed, and the uncertainty in the amount of recoverable oil in a reservoir can be resolved in two different ways, namely drilling a sufficient number of wells or using the reservoir for production for long enough.

\subsection{Incorporating conditionally observed information in rooted junction trees}
\label{subsection:cnac_obs}
The conditional arcs are designed to describe conditionally observed information in influence diagrams. However, they are a general representation of the concept, not a modeling solution. In what follows, we present two alternative approaches for incorporating this concept into the RJT models, which ultimately enables solving Type 3 endogenously uncertain stochastic problems. The first approach employs \emph{observation nodes}. These have been used in influence diagram literature, in problems such as the used car buyer problem \citep{howard2005influence} or the oil wildcatter problem \citep{raiffa1968decision}. The second approach utilizes conditional non-anticipativity constraints (C-NACs), which are used in stochastic programming for modeling decision-dependent information structures. While C-NACs have been somewhat widely used in scenario tree-based stochastic programming, as discussed in Section \ref{section:uncertainty}, they have not been applied in the context of influence diagram models and Type 3 endogenous uncertainty. To fill this gap, we show that conditional non-anticipativity constraints can be incorporated in the mixed-integer reformulation \eqref{eq:rjt-obj}-\eqref{eq:rjt-vars} for modeling and solving limited memory influence diagrams. In Section \ref{section:experiments}, the computational performance of C-NACs and observation nodes is compared.

\subsubsection{Observation nodes}

Observation nodes portray how the decision maker observes the information. By enforcing that earlier decisions affect the probability distribution of the observations, Type 2 uncertainty is effectively transformed into Type 1 uncertainty, making it directly amenable to influence diagrams and their associated analysis techniques. This approach is also used in \citet{Salo2022}. While observation nodes have been previously used in the literature on influence diagrams, the specific connection between conditionally observed information and observation nodes has not been explicitly discussed nor has the computational impact of incorporating conditional information been assessed.

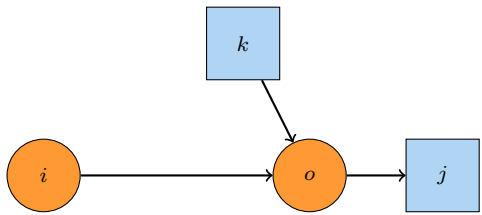
\begin{figure}[H]
\centering
\begin{tikzpicture}
    [decision/.style={fill=blue!40, draw, minimum size=2.5em, inner sep=2pt}, 
    chance/.style={circle, fill=orange!80, draw, minimum size=2.5em, inner sep=2pt},
    value/.style={diamond, fill=teal!80, draw, minimum size=2.5em, inner sep=2pt},
    optimization/.style={ellipse, fill=teal!80, draw, minimum size=2em, inner sep=2pt},
    scale=1.75, font=\scriptsize]
    
    \node[chance]   (L) at (-0.25, 1.5)  {$I_{pre}$};
    \node[chance]   (L1) at (0.75, 2.25)  {$R_1$};
    \node[decision] (1) at (2.25, 2.25)  {$A_1$};
    \node[chance]   (F) at (3.25, 1.5)  {$I_{post}$};
    \node[value]    (T) at (4.25, 1.5)  {$U_{int}$};   
    \node[value]    (V1) at (3.25, 2.5)  {$V_1$};      
    \node    (2) at (1.5, 0.75)  {...};    
    \node[decision]  (I) at (1.5, 3.5) {$D_1$};
    \node[chance] (O) at (1.5,2.25) {$O_1$};
    \draw[->, thick] (L1) -- (O);
    \draw[->, thick] (I) --  (O);
    \draw[->, thick] (O) --  (1);
    \draw[->, thick] (L) -- (L1);
    \draw[->, thick] (L) -- (F);
    \draw[->, thick] (1) -- (F);
    \draw[->, thick] (2) -- (F);
    \draw[->, thick] (F) -- (T);
    \draw[->, thick] (1) -- (V1);
    \draw[->, thick] (L) -- (2);
    \draw[->, thick] (L) to[out=60,in=135,distance=1.5cm] (V1);
\end{tikzpicture}
\caption{Replacing a distinguishability arc with an observation node in the example from Fig. \ref{fig:dist_arc_basic}.}
\label{fig:dist_arc_obs}
\end{figure}

In effect, each conditional arc is replaced with an observation node, as illustrated in Fig. \ref{fig:dist_arc_obs}. The information set of the observation node $O_1$ is the union of the chance node $R_1$ and the distinguishability set $T_{R_1,A_1} = \{D_1\}$, and the state space is $S_{R_1} \cup \text{``no observation"}$. Then, the observation node replaces the node $R_1$ in the information set of $A_1$, controlling whether or not the information in $R_1$ is available when making the decision in $A_1$. A benefit of this approach is that the modifications are done to the influence diagram, and the RJT formed from the resulting diagram can immediately be used in \eqref{eq:rjt-obj}-\eqref{eq:rjt-vars}.

One can also utilize the ideas of \cite{Ruszczynski1997} and \cite{apap2017} for modeling conditional information revelation, in which the information structure can be connected to the decisions by using disjunctive conditional non-anticipativity constraints (C-NACs). These constraints are similar to the more traditional non-anticipativity constraints \citep{rockafellar1991scenarios} in stochastic programming, but the constraints are only imposed if the distinguishability conditions between pairs of scenarios are not satisfied. These constraints are designed for scenario tree-based MIP formulations of decision problems and, hence, have not been directly amenable to influence diagram models thus far.

\subsubsection{Conditional nonanticipativity constraints (C-NAC)}

To integrate C-NACs into our model, we first ease the notation by supplementing the conditional arcs with the \emph{conditional information set} $I_c(j) = \{i \in N^C \cup N^D \mid (i,j) \in A_c \}$ to represent the conditionally available information at node $j \in N^D$. In practice, the conditional information set is a part of the information set of $j$ ($I_c(j) \subseteq I(j)$) and the C-NACs then control whether or not this information is available when making the decision represented by node $j$. Note that to formulate the C-NACs, the related distinguishability sets $T_{i,j}$ for all $i \in I_c(j)$ must also be contained in the cluster $C_j$ as the constraints will be written for a single cluster and propagated across the RJT via constraints \eqref{eq:rjt-prob}-\eqref{eq:rjt-dec}. Since any algorithm converting an influence diagram to an RJT should ensure the property $I(j) \in C_j$, this can be achieved by treating them as a part of the information set of $j$ when converting the influence diagram to an RJT. The MIP \eqref{eq:rjt-obj}-\eqref{eq:rjt-vars} can then be formulated from the RJT as described in Section \ref{subsec:rjt}, and the conditionally observed information is modeled by adding constraints \eqref{eq:C-NAC} to the resulting formulation.

C-NACs are used to enforce conditional non-anticipativity when the conditional information states $s_{I_c(j)}$ differ between cluster states $s_{C_j}$ and $s'_{C_j}$. As discussed earlier, if the conditional information states differ (i.e., node $i \in I_c(j)$ has different states in $s_{C_j}$ and $s'_{C_j}$), distinguishability is dependent on the corresponding condition(s) $F_{i,j}(s_{T_{i,j}})$. This distinguishability of two cluster states at node $j \in N^D$ can be formulated as a Boolean variable $f_j^{s_{C_j},s'_{C_j}}$, defined as 
$$f_j^{s_{C_j},s'_{C_j}}=\begin{cases}
    True & \ \text {if } \exists i \in I_c(j), s_i \neq s'_i, F_{i,j}(s_{T_{i,j}})=1 \\
    False & \ otherwise.
\end{cases}$$
The value of $f_j^{s_{C_j},s'_{C_j}}$ is \textit{True} (i.e., 1) if any node $i$, for which $s_i \neq s'_i$, is observed when making the decision in $j$ and thus makes scenarios $s_{C_j}$ and $s'_{C_j}$ distinguishable at node $j$, and \textit{False} (i.e., 0) otherwise. 

Finally, we extend the definition of the local decision strategy $Z_j(s_{I(j)})$ and the corresponding binary variables $z(s_j \mid s_{I(j)})$ to differentiate between the conditional information set $I_c(j)$ and the ``non-conditional" information set $I_n(j)$ containing the nodes that are always observed at $j$. If the value of $f_j^{s_{C_j},s'_{C_j}}$ is \textit{False}, a constraint is required to force the local strategies $Z_j(s_{I(j)},s_{I_c(j)})$ and $Z_j(s_{I(j)},s'_{I_c(j)})$ to be the same. Note that we only consider pairs of strategies for which $s_{I_n(j)}$ is the same, as different non-conditional information states are always distinguishable, precluding the need for C-NACs. Finally, we can define C-NACs in the context of our model as 
\begin{align}
  \neg f_j^{s_{C_j},s'_{C_j}} \Longrightarrow z(s_j \mid s_{I_n(j)},s_{I_c(j)}) = z(s_j \mid s_{I_n(j)},s'_{I_c(j)}), &\, \forall j \in N^D, s_j \in S_j,\\ & s_{I_n(j)} \in S_{I_n(j)}, s_{I_c(j)}, s'_{I_c(j)} \in S_{I_c(j)} \nonumber.
\end{align}

In light of the above, the C-NACs for binary variables $z$ can also be conveniently written as
\begin{align}\label{eq:C-NAC}
  | z(s_j \mid s_{I_n(j)},s_{I_c(j)}) - z(s_j \mid s_{I_n(j)},s'_{I_c(j)}) | \le f_j^{s_{C_j},s'_{C_j}}, &\, \forall j \in N^D, s_j \in S_j,\\ & s_{I_n(j)} \in S_{I_n(j)}, s_{I_c(j)}, s'_{I_c(j)} \in S_{I_c(j)} \nonumber.
\end{align}
Notice that the absolute value function used in the left-hand side of \eqref{eq:C-NAC} can be trivially linearized without significantly increasing the model complexity. This constraint states that for each decision node $j \in N^D$, if the non-conditional information states $s_{I_n(j)}$ are the same for two cluster states $s_{C_j}$ and $s'_{C_j}$, and conditionally revealed information does not make these cluster states distinguishable either, the corresponding local decision strategies represented by the $z$-variables must be the same.

In practice, the main challenge with using C-NACs is that the number of constraints \eqref{eq:C-NAC} quickly becomes overwhelmingly large. With this in mind, \citet{apap2017} present several properties that C-NACs possess that can be exploited to reduce the number of such constraints. By making use of these C-NAC reduction properties, representing the decision-dependent information structure within the RJT model is likely to be more compact with C-NACs than the corresponding model using observation nodes, as C-NACs do not require a separate node for modeling the conditional observation.

One limitation of the C-NACs is their implicit assumption that a decision reveals the conditionally observed information with certainty. As such, they cannot be directly employed if the decision reveals the information with a probability of less than one, unlike observation nodes. However, this can easily be circumvented by adding a chance node representing the success/failure of the observation decision and using this new node in the distinguishability function $f$ in \eqref{eq:C-NAC}.

Considering the MIP model \eqref{eq:rjt-obj}-\eqref{eq:rjt-vars}, observation nodes result in the conditional probability distribution associated with the observation node setting some probabilities $\mu$ to zero. For example, the probability of observing a state other than ``no observation'' is zero if the distinguishability condition is not fulfilled. On the other hand, C-NACs result in fewer clusters and do not use such extended conditional probability distributions to set $\mu$-variables to zero. Instead, they operate by only imposing equality constraints between decision strategies with indistinguishable information, typically yielding a smaller model. This hypothesis, and the computational performance of the two approaches, are explored in the computational experiments in Section \ref{section:experiments}.

\subsection{Subdiagram decomposition}
\label{section:decomp}

While the MIP formulations in \citet{parmentier2020integer} and \citet{Salo2022} offer excellent modeling flexibility for influence diagrams with additional constraints and alternative objective functions, they can only accommodate problems where all decisions have a discrete and finite set of alternatives. However, many problem settings involve decisions that are more naturally modeled as continuous variables. For instance, the climate change cost-benefit analysis in Section \ref{section:score} that motivated our developments considers the optimal emission levels in the future. A straightforward approach to make the modeling of these decisions amenable to the proposed MIP formulations is to discretize them with a potentially large number of states. The two main disadvantages of this approach are suboptimal solutions if the discretization fails to accurately represent the decision space, and the increase in model size. In this subsection, we first present the conditions for applying the decomposition using the running N-monitoring example, before discussing the practical aspects of using the decomposition.

To allow for continuous decisions, we build on the ideas in \citet{lee2021submodel}, where the authors present a submodel decomposition for influence diagrams based on partitioning the diagram into subdiagrams whose optimal strategies are independent of each other. However, they only consider influence diagrams with discrete decisions, similarly to \citet{Salo2022} and \citet{parmentier2020integer}. Our key contribution to exploiting this decomposition is that if we can partition the original diagram into subdiagrams in a way that decision-dependent probabilities and continuous decision variables do not appear in the same subdiagram, we can model problems with both Type 3 endogenous uncertainty and continuous decision variables. For that, we must be able to decompose the diagram into subdiagrams with decision-dependent probabilities and discrete decisions and subdiagrams with continuous decisions but no decision-dependent probabilities. If that is possible, the former can be solved using the MIP formulations presented in this paper, designed for Type 3 endogenously uncertain problems. The latter can be handled using standard stochastic programming modeling methods that can model exogenous uncertainty.

While it may seem somewhat restrictive to not allow decision-dependent probabilities and continuous decisions in the same subproblem, this structure is in line with the examples in \citet{hellemo2018}, where first-stage decisions affect the probability distributions in later stages. Additionally, many typical stochastic problems, such as facility location \citep{baron2008facility} and unit commitment \citep{zheng2014stochastic} involve a structure where the discrete decisions in earlier stages affect the continuous decisions in later stages. 

We illustrate the decomposition with an extended version of the N-monitoring problem, presented in Fig. \ref{fig:nmon-models}. Let us assume that instead of the previous setting where nodes $A_i$ represented the decision of whether to sell the product at location $i$, the regional representatives now decide their production level $O_i$, with $A_i$ being a discrete decision determining the upper and lower bound for the continuous variable represented by $O_i$. The intuition behind this is that instead of only considering the impact of introducing the product in each area, the company has conducted a more careful analysis on the relationship between production amounts and the post-decision interest level $I_{post}$, assuming that large product launches have a greater impact on overall interest than producing few units that end up in only a handful of stores. However, instead of trying to define separate conditional probabilities in $I_{post}$ for producing any arbitrary number of products in each area, the company uses a discretization represented by nodes $A_i$. These discrete decisions, e.g., ``500-1000 units produced in area 1", are then used to estimate the effect on the overall public interest level $I_{post}$, for which we give a utility value determined by $U_{int}$. The exact production amounts represented by $O_i$ then affect the regional profits $V_i$. In Fig. \ref{fig:nmon-bad}, the continuous decisions $O_i$ have an effect on the utility $U_{int}$; in Fig. \ref{fig:nmon-ext} they do not. This difference will be used to discuss the limitations of the proposed decomposition approach.

\begin{figure}[!ht]
\centering
     \begin{subfigure}[b]{0.48\textwidth}
        \centering
        \resizebox{\textwidth}{!}{%
        \begin{tikzpicture}
            [decision/.style={fill=blue!40, draw, minimum size=2.5em, inner sep=2pt}, 
            chance/.style={circle, fill=orange!80, draw, minimum size=2.5em, inner sep=2pt},
            value/.style={diamond, fill=teal!80, draw, minimum size=2.5em, inner sep=2pt},
            optimization/.style={ellipse, fill=teal!80, draw, minimum size=2em, inner sep=2pt},
            scale=1.5, font=\scriptsize]
            
            \node[chance]   (L) at (-0.25, 1.5)  {$I_{pre}$};
            \node[chance]   (L1) at (0.75, 2.25)  {$R_1$};
            \node[chance]   (L2) at (0.75, 0.75)  {$R_2$};
            \node[decision] (1) at (1.75, 2.25)  {$A_1$};
            \node[decision] (2) at (1.75, 0.75)  {$A_2$};
            \node[decision] (O1) at (2, 3)  {$O_1$};
            \node[decision] (O2) at (2, 0)  {$O_2$};
            \node[chance]   (F) at (2.75, 1.5)  {$I_{post}$};
            \node[value]    (T) at (3.75, 1.5)  {$U_{int}$};   
            \node[value]    (V1) at (3, 3)  {$V_1$};   
            \node[value]    (V2) at (3, 0)  {$V_2$};     
            \draw[->, thick] (L) -- (L1);
            \draw[->, thick] (L) -- (L2);
            \draw[->, thick] (L) -- (F);
            \draw[->, thick] (L1) -- (1);
            \draw[->, thick] (L2) -- (2);
            \draw[->, thick] (1) -- (F);
            \draw[->, thick] (2) -- (F);
            \draw[->, thick] (F) -- (T);
            \draw[->, thick] (L1) -- (O1);
            \draw[->, thick] (L2) -- (O2);
            \draw[->, thick] (1) -- (O1);
            \draw[->, thick] (2) -- (O2);
            \draw[->, thick] (O1) -- (V1);
            \draw[->, thick] (O2) -- (V2);
            \draw[->, thick] (L) to[out=60,in=135,distance=1.5cm] (V1);
            \draw[->, thick] (L) to[out=-60,in=-135,distance=1.5cm] (V2);
        \end{tikzpicture}
        }%
        \caption{N-monitoring problem with continuous decisions}
        \label{fig:nmon-ext}
     \end{subfigure}
     \hfill
     \begin{subfigure}[b]{0.48\textwidth}
        \centering
        \resizebox{\textwidth}{!}{%
        \begin{tikzpicture}
            [decision/.style={fill=blue!40, draw, minimum size=2.5em, inner sep=2pt}, 
            chance/.style={circle, fill=orange!80, draw, minimum size=2.5em, inner sep=2pt},
            value/.style={diamond, fill=teal!80, draw, minimum size=2.5em, inner sep=2pt},
            optimization/.style={ellipse, fill=teal!80, draw, minimum size=2em, inner sep=2pt},
            scale=1.5, font=\scriptsize]
            
            \node[chance]   (L) at (-0.25, 1.5)  {$I_{pre}$};
            \node[chance]   (L1) at (0.75, 2.25)  {$R_1$};
            \node[chance]   (L2) at (0.75, 0.75)  {$R_2$};
            \node[decision] (1) at (1.75, 2.25)  {$A_1$};
            \node[decision] (2) at (1.75, 0.75)  {$A_2$};
            \node[decision] (O1) at (2, 3)  {$O_1$};
            \node[decision] (O2) at (2, 0)  {$O_2$};
            \node[chance]   (F) at (2.75, 1.5)  {$I_{post}$};
            \node[value]    (T) at (3.75, 1.5)  {$U_{int}$};   
            \node[value]    (V1) at (3, 3)  {$V_1$};   
            \node[value]    (V2) at (3, 0)  {$V_2$};     
            \draw[->, thick] (L) -- (L1);
            \draw[->, thick] (L) -- (L2);
            \draw[->, thick] (L) -- (F);
            \draw[->, thick] (L1) -- (1);
            \draw[->, thick] (L2) -- (2);
            \draw[->, thick] (1) -- (F);
            \draw[->, thick] (2) -- (F);
            \draw[->, thick] (F) -- (T);
            \draw[->, thick] (L1) -- (O1);
            \draw[->, thick] (L2) -- (O2);
            \draw[->, thick] (1) -- (O1);
            \draw[->, thick] (2) -- (O2);
            \draw[->, thick] (O1) -- (V1);
            \draw[->, thick] (O2) -- (V2);
            \draw[->, thick] (O1) -- (T);
            \draw[->, thick] (O2) -- (T);
            \draw[->, thick] (L) to[out=60,in=135,distance=1.5cm] (V1);
            \draw[->, thick] (L) to[out=-60,in=-135,distance=1.5cm] (V2);
        \end{tikzpicture}
        }%
        \caption{N-monitoring problem with continuous decisions and added influence on $U_{int}$}
        \label{fig:nmon-bad}
     \end{subfigure}
     \caption{Two different influence diagram representations for the N-monitoring problem.}
     \label{fig:nmon-models}
\end{figure}

The submodel decomposition can be seen as taking a part of the original diagram and converting it into a value node representing the submodel. However, one cannot simply choose a subset of nodes and make the corresponding part of the diagram a subproblem. Instead, the subproblem must be \emph{stable}, such that the set of relevant hidden variables does not contain any unobserved decision variable \citep{lee2021submodel}. The key observation from \citet{lee2021submodel} is that for a stable submodel, the task of maximizing expected utility within the submodel is independent of the decision strategies outside the submodel.

More formally, a relevant hidden node for decision node $d$ is a node $a$ such that the descendant value nodes $N^V_d$ of $d$ (that is, value nodes to which there is a directed path from $d$) are not d-separated \citep{pearl2009causality} from $a$ by either $d$ or $I(d)$. The lack of separation means that $a$ influences the same value node as $d$, and the influence cannot be determined only knowing the states of $I(d)$ and $d$. We illustrate this in Fig. \ref{fig:nmon-bad} with a slightly modified version of the N-monitoring problem that does not lend itself to decomposition. A submodel consisting of nodes $O_1$, $V_1$ and $U_{int}$ is not stable in Fig. \ref{fig:nmon-bad}, because the path $O_2 \rightarrow U_{int}$ is not blocked by $O_1$ or $I(O_1)=\{R_1, A_1\}$ and knowing the decision strategy in $O_2$ is thus relevant when making the decision in $O_1$. For a more thorough treatment of the topic, we refer the reader to \citet{lauritzen2001representing} and \citet{lee2021submodel}.

In Fig. \ref{fig:nmon-ext}, a suitable decomposition is to separate the continuous decision $O_i$ and the associated value node $V_i$ into a subproblem for each $i \in \{1,...N\}$. The utility value of these subproblems depends on $I_{pre}$, $R_i$ and $A_i$, and the decisions $A_j$ and $O_j$ for $j \neq i$ are not relevant in the subproblem $i$. However, since the state of $I_{pre}$ is not known when making the decision $O_i$, we solve stochastic subproblems $SP_i(s_{R_i},s_{A_i})$ parametrized by the states of the nodes in $I(O_i)$. From these subproblems, optimal decision strategies for nodes $O_i$, independent of $s_{I_{pre}}$ can be obtained. These subproblems consist of a continuous decision on how much to produce, and the unknown demand represented by node $I_{pre}$ affects the resulting profit, making each subproblem a stochastic newsvendor problem \citep{petruzzi1999pricing} where the production level is constrained by the decision $s_{A_i}$ and an imperfect estimate of demand is given by $s_{R_i}$. 

Naturally, the structure of the subproblem affects the computational complexity of solving the full problem, and the more complex (e.g., nonlinear and/or nonconvex) the subproblems, the more technically and computationally challenging it is to solve them to optimality. In theory, the only requirement for the subproblem is that it can be solved, as only the optimal subproblem solutions are needed in the main problem, different from decomposition approaches such as Benders \citep{benders1962partitioning} where strong duality is required to hold for the subproblem. 

While this decomposition can be a powerful tool for solving influence diagrams, one of the main disadvantages in the context of our MILP models is that, as most solution methods for IDs, it generally only works for ``pure" expected utility maximization. For example, if we wanted to impose a chance constraint on the sum of utilities in nodes $V_i$ and $U_{int}$ in Fig. \ref{fig:nmon-ext}, all of the local decision strategies would have to be considered together, making a decomposition approach infeasible. On the other hand, conditionally observed information can be straightforwardly handled if the node and its conditional information set are in the same subdiagram. If not, conditional non-anticipativity would have to be enforced between subproblems by combining the subproblems into a larger, stochastic subproblem. For instance, if the reports $R_i$ in Fig. \ref{fig:nmon-ext} were conditionally observed, the subproblems $SP_i(s_{I_{pre}},s_{R_i},s_{A_i})$ for all $s_{R_i} \in S_{R_i}$ would need to be combined into one subproblem to enforce C-NACs for the decision in $O_i$.

\section{Computational experiments}
\label{section:experiments}

All problems are solved using 8 threads on an Intel E5-2680 CPU at 2.5GHz and 16GB of RAM, provided by the Aalto University School of Science ``Science-IT'' project. The code was implemented in Julia v1.10.3 \citep{Julia-2017} using Gurobi v11.0.0 \citep{gurobi} and JuMP v1.22.1 \citep{DunningHuchetteLubin2017}. All the code used in the computational experiments is available at a GitHub repository \citep{repo}.

\begin{figure}[!ht]
\centering
     \begin{subfigure}[b]{0.8\textwidth}
        \centering
        \includegraphics[width=\textwidth]{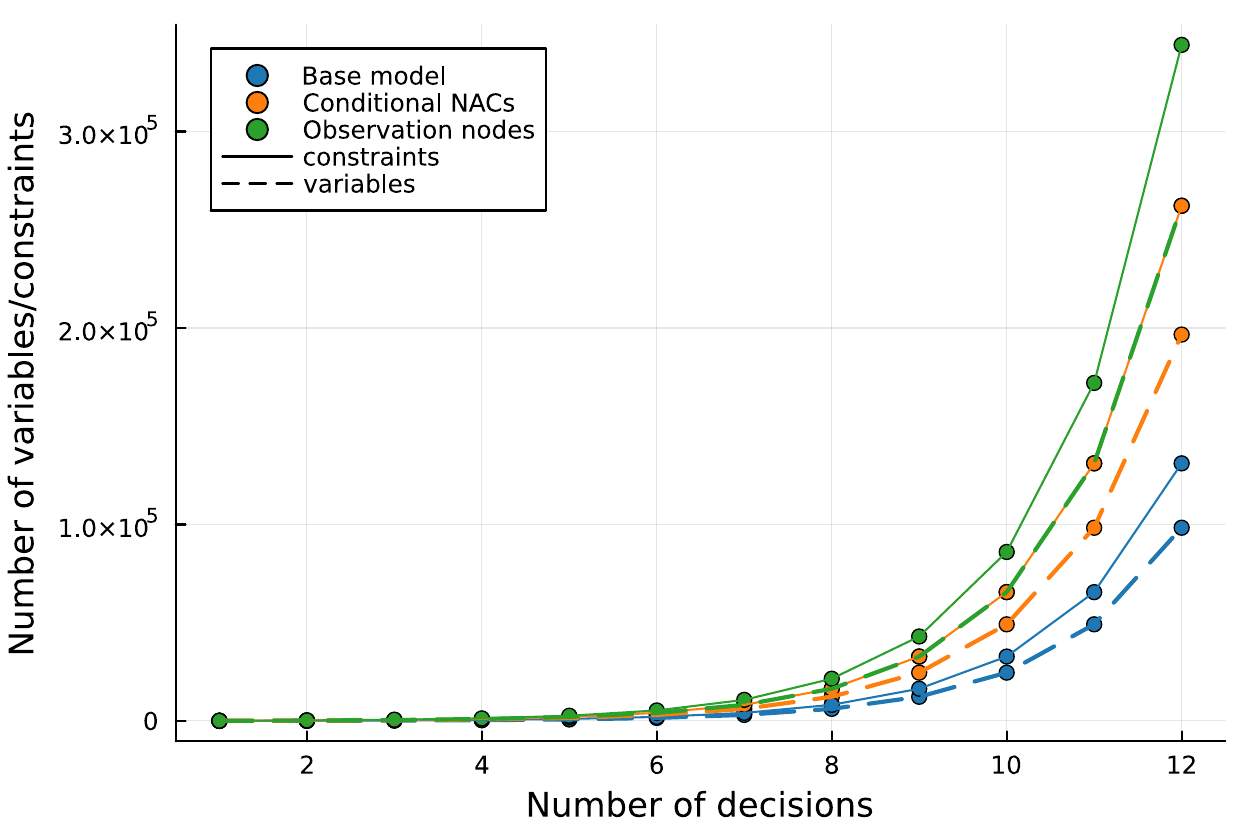}
        \caption{Model sizes}
        \label{fig:nmon-modelsizes}
     \end{subfigure}
     \begin{subfigure}[b]{0.8\textwidth}
        \centering
        \includegraphics[width=\textwidth]{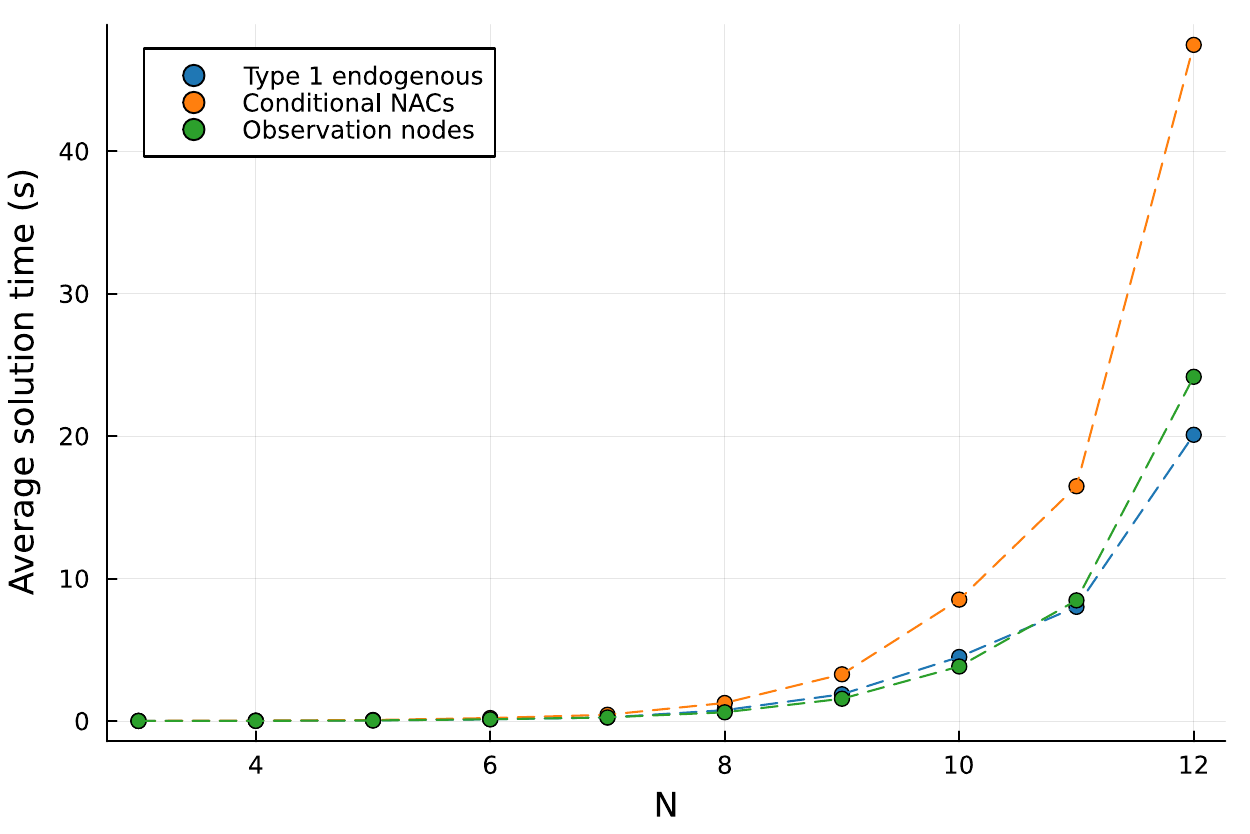}
        \caption{Solution times}
        \label{fig:nmon-soltimes}
     \end{subfigure}
     \caption{Model sizes and solution times for the N-monitoring problem with conditionally observed information using C-NACs and observation nodes. ``Base model'' refers to the model with no conditional observation (Fig. \ref{fig:nmon}).}
     \label{fig:nmon-experiments}
\end{figure}

First, we solve the N-monitoring problem with conditionally observed information, varying the number of decision makers. The model sizes and average solution times over 50 randomly parametrized instances are compared between models using observation nodes and C-NACs, and against the original N-monitoring problem (Fig. \ref{fig:nmon}). The randomized instance generation for the conditional N-monitoring problem is described in Appendix \ref{sec:app_nmon}, and the results are presented in Fig. \ref{fig:nmon-experiments}. 

We can see that the number of constraints and variables grows exponentially as more decision makers are considered, approximately doubling as $N$ is increased by one. As expected, the model with C-NACs is significantly smaller than the model using observation nodes when modeling conditionally observed information. However, the solution times for observation nodes are close to those in the model with no conditional information revelation, while the smaller C-NAC model results in significantly longer solution times. Further experimentation shows that finding an initial feasible solution takes significantly longer for the C-NAC model. This is likely related to the structure of the RJT model \eqref{eq:rjt-obj}-\eqref{eq:rjt-vars} in that the zero probabilities in constraint \eqref{eq:rjt-prob} resulting from observation nodes are efficiently propagated across clusters, while C-NACs, being additional constraints, make the structure of the model more difficult for the solver. This is in line with the observation that the RJT model \eqref{eq:rjt-obj}-\eqref{eq:rjt-vars} by itself provides a strong linear relaxation \citep{parmentier2020integer,hankimaa2023}.

In order to explore the effect of the decomposition and continuous decision variables, we use the version of the N-monitoring problem described in Section \ref{section:contributions} with conditionally observed information modeled using C-NACs. We solve 50 instances of the 4-monitoring problem with $k$ randomly selected discrete decision alternatives representing the continuous decision and compare the results to those obtained from the continuous version of the model in Table \ref{tbl:nmon-results}. As discussed earlier in Section \ref{section:decomp}, the subproblems are instances of the newsvendor problem. For such problems, under mild assumptions (positive profit from selling products within the realized demand, negative profit for products exceeding demand), production levels greater than maximum demand or lower than minimum demand are always suboptimal. Thus, we sample the possible decisions (i.e., the states for the decision nodes being discretized) from a uniform distribution between the minimum and maximum demand.

\begin{table}[!ht]
\centering
\begin{tabular}{l|lllllllllll}
\# of discrete alternatives $k$ & 1    & 2    & 3    & 4     & 6    & 8    & 10   & $\infty$ \\ \hline
solution time (s)               & 0.03 & 0.12 & 0.36 & 1.56 & 5.73 & 24.57 & 69.10 & 0.03 \\
relative objective value        & 0.908 & 0.949 & 0.964 & 0.963 & 0.975 & 0.987 & 0.997 & 1   
\end{tabular}
\caption{Average solution times and objective values for discrete and continuous versions of the problem. Objective values are scaled such that the actual optimal value is 1. The decomposed model with continuous decision is denoted as having infinitely many decision alternatives.}
\label{tbl:nmon-results}
\end{table}

It can be seen from the results in Table \ref{tbl:nmon-results} that a finer discretization results in solutions that are closer to the optimal objective value, but with an increasing computational cost. This rapid increase in solution times requires a tradeoff between computational tractability and solution quality. Using the subdiagram decomposition to separate the continuous decision from the diagram results in a solution time similar to the coarsest discretizations and the best possible objective value. This shows that while discretizing continuous variables might be a straightforward solution, it suffers from both high solution times and suboptimal solutions. If the modeler has prior information on the nature of the optimal decisions, the discretization could be more efficient than random sampling from a uniform distribution, but the same issues would persist, albeit lessened.

\section{Cost-benefit analysis for climate change mitigation}
\label{section:score}

\subsection{Model description}

To illustrate the setting of Type 3 endogenous uncertainty with continuous decision variables, we now consider the cost-benefit analysis on mitigating climate change under uncertainty \citep[see, e.g., ][]{Ekholm2018}. Climate change is driven by greenhouse gas (GHG) emissions and can be mitigated by reducing these emissions, which incurs costs. Still, mitigation reduces the negative impacts of climate change, referred to as climate damage. In cost-benefit analysis, the objective is to minimize the discounted sum of mitigation costs and climate damage over a long time horizon. However, multiple uncertainties complicate the analysis.

Here, we consider three salient uncertainties involving both decision-dependent probabilities (Type 1) and conditionally observed information (Type 2). Moreover, some decision nodes involve continuous variables. The resulting problem is a multi-stage mixed-integer nonlinear problem (MINLP) with Type 3 endogenous uncertainty, thus demonstrating the proposed novel features to the rooted junction tree framework described in this paper.

For the mitigation costs, the model considers that technological R\&D can be conducted to decrease the costs of bioenergy with carbon capture and storage (BECCS). These R\&D decisions can take place at three intensity levels over two distinct stages. The first stage is a choice between low or medium R\&D effort. The low-effort choice represents a business-as-usual perspective, which carries throughout the decision process. If the medium effort is chosen, one observes whether the R\&D looks promising or not, and can then decide whether to continue with the medium or switch to a higher R\&D effort. The three R\&D effort levels and whether the development seems promising or not all affect the probabilities for achieving either low, medium or high mitigation costs later during the century. 

The presented model is an extension from \citet{EkholmBaker2021}, which in turn is a simplification from the SCORE model \citep{Ekholm2018}. Compared to the formulation proposed here, these earlier analyses have assumed that the uncertainties regarding climate parameters are resolved exogenously over time and the uncertainty in mitigation cost is modeled through separate scenarios. The details of the model structure and parametrization are summarized here and described in detail in Appendix \ref{sec:app_score}. It should be noted that the model presented here is an oversimplified representation of reality, and obtaining suitable, defendable estimates for its parameters is an extremely challenging task. Thus, this case study should be viewed as illustrative. Nevertheless, the values used in this paper are based on expert elicitation and the results are in line with existing literature, suggesting that the proposed framework can be used for problems comprising different types of endogenous uncertainty.

\begin{figure}[!ht]
\centering
\resizebox{\textwidth}{!}{%
\begin{tikzpicture}
    [decision/.style={fill=blue!40, draw, minimum size=2.5em, inner sep=2pt}, 
    chance/.style={circle, fill=orange!80, draw, minimum size=2.5em, inner sep=2pt},
    value/.style={diamond, fill=teal!80, draw, minimum size=2.5em, inner sep=2pt},
    optimization/.style={ellipse, fill=teal!80, draw, minimum size=2em, inner sep=2pt},
    scale=1.625, font=\scriptsize]
     \node[decision, align=center] (D1) at (0, 2)   {$D_{Dmg}$\\9};
     \node[decision, align=center] (D2) at (0, 1)   {$D_{T1}$\\1};
     \node[decision, align=center] (D3) at (0, 0)   {$D_{CS}$\\6};
     \node[decision, align=center] (D4) at (2.5, 1.5)   {$D_{E1}$\\5};
     \node[decision, align=center] (D5) at (2.5, 0.5)   {$D_{T2}$\\3};
     \node[decision, align=center] (D6) at (4.5, 1)   {$D_{E2}$\\12};
     \node[decision, align=center] (D7) at (6.0, 1)   {$D_{E3}$\\13};
     \node[chance, align=center]   (C1) at (3.25, 2)   {$O_{Dmg}$\\10};
     \node[chance, align=center]   (C2) at (1, 1)   {$O_{T1}$\\2};
     \node[chance, align=center]   (C3) at (3.25, 0)   {$O_{CS}$\\7};
     \node[chance, align=center]   (C4) at (3.5, 1)   {$O_{T2}$\\4};
     \node[chance, align=center]   (O1) at (4, 3)   {$R_{DMG}$\\11};
     \node[chance, align=center]   (O2) at (4, -1)   {$R_{CS}$\\8};
     \draw[] (4.1,2.1) pic{dist_arc={int1, , ,rotate=105}};
     \draw[] (4.1,-0.1) pic{dist_arc={int2, , ,rotate=75}};
     \node[value, align=center] (SC) at (7.5, 1)   {Costs, \\damages};
     \draw[->, thick] (O1) -- (D7);
     \draw[->, thick] (O1) -- (SC);
     \draw[->, thick] (O2) -- (D7);
     \draw[->, thick] (O2) -- (SC);
     \draw[->, thick] (D1) -- (C1);
     \draw[->, thick] (D2) -- (C2);
     \draw[->, thick] (D3) -- (C3);
     \draw[->, thick] (C2) -- (D4);
     \draw[->, thick] (C2) -- (D5);
     \draw[->, thick] (C2) -- (C4);
     \draw[->, thick] (D5) -- (C4);
     \draw[->, thick] (C4) -- (D6);
     \draw[->, thick] (C4) to[out=45,in=135,distance=0.3cm] (D7);
     \draw[->, thick] (D4) to[out=0,in=120,distance=0.5cm] (SC);
     \draw[->, thick] (D6) -- (D7);
     \draw[->, thick] (D6) to[out=45,in=135,distance=0.3cm] (SC);
     \draw[->, thick] (D7) -- (SC);
     
     \draw[thick] (C1) -- (int1.195);
     \draw[thick] (O1) -- (int1.105);
     \draw[->, thick] (int1.285) -- (D6);
     
     \draw[thick] (C3) -- (int2.165);
     \draw[thick] (O2) -- (int2.255);
     \draw[->, thick] (int2.75) -- (D6);
     
     \draw[->, thick] (0,-2) -- (8,-2);
     \draw[-, thick] (0.0,-1.9) -- (0.0,-2.1);
     \draw[-, thick] (2.5,-1.9) -- (2.5,-2.1);
     \draw[-, thick] (4.5,-1.9) -- (4.5,-2.1);
     \draw[-, thick] (6.0,-1.9) -- (6.0,-2.1);
     \draw[-, thick] (7.5,-1.9) -- (7.5,-2.1);
     \node[] at (0.0,-2.2) {2020};
     \node[] at (2.5,-2.2) {2030};
     \node[] at (4.5,-2.2) {2050};
     \node[] at (6.0,-2.2) {2070};
     \node[] at (7.5,-2.2) {2100};
\end{tikzpicture}
}%
\caption{Influence diagram of the climate change cost-benefit problem with endogenous uncertainty due to R\&D. We assume that all prior decisions and uncertainty realizations apart from the conditionally observed parameters are remembered when making decisions, but omit the arcs for clarity. Additionally, the individual value nodes associated with decision nodes 1, 3, 6 and 9 are omitted for clarity. The nodes are numbered according to the topological order used for forming the RJT.}
\label{fig:score}
\end{figure}
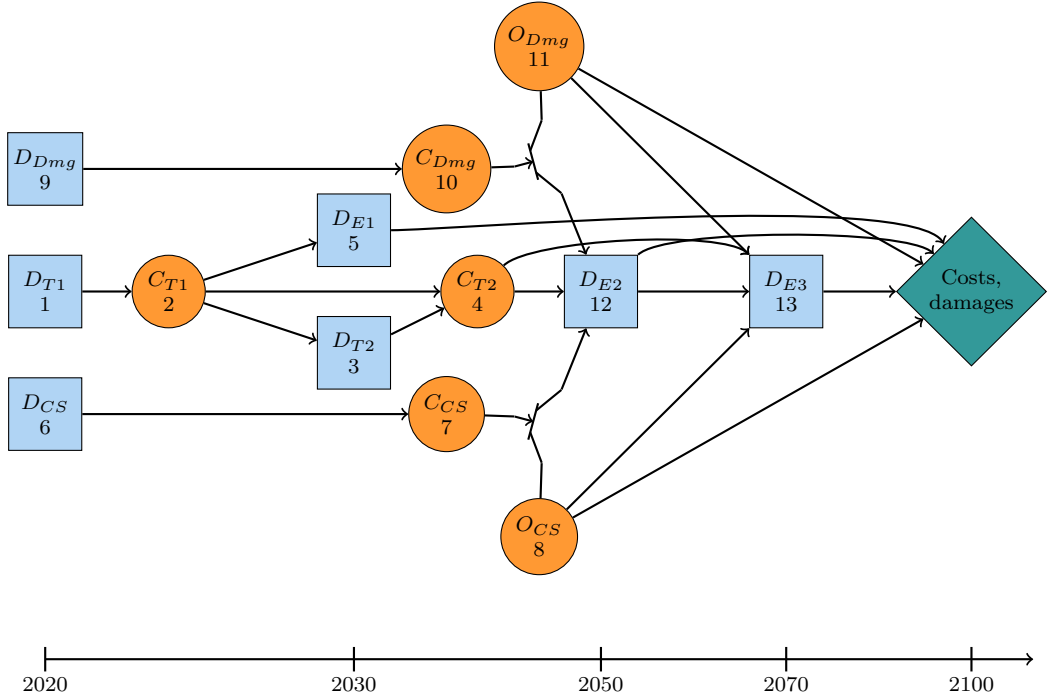

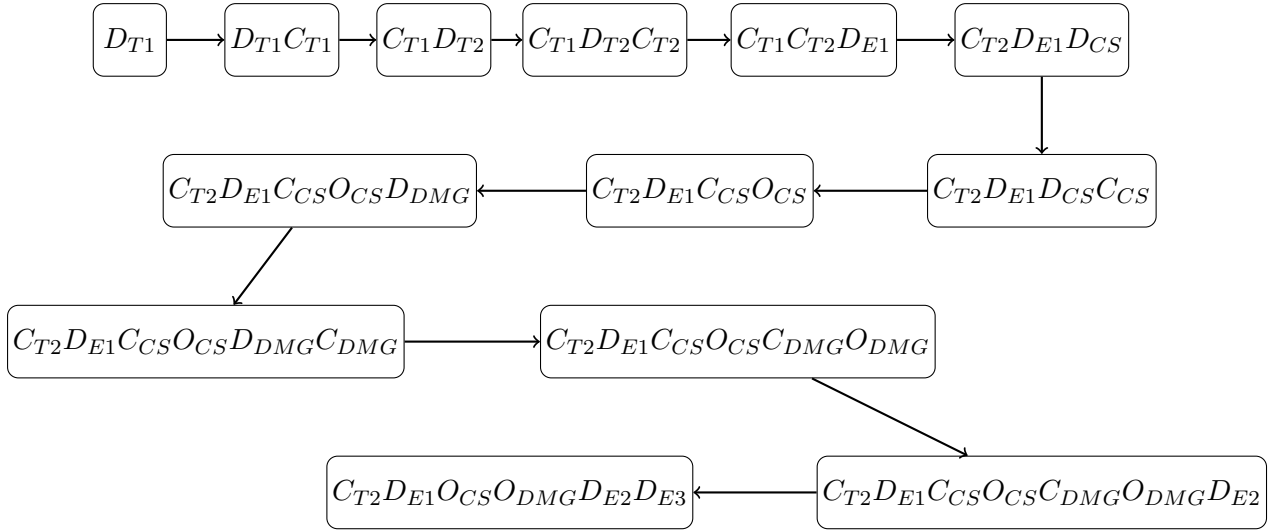
\begin{figure}[!ht]
\centering
\resizebox{\textwidth}{!}{%
\begin{tikzpicture}
    [cluster/.style={fill=white!80, draw, minimum size=2.5em, inner sep=2pt, rounded corners}]
    \node[cluster] (1) at (0, 0)   {$D_{T1}$};
    \node[cluster] (2) at (2, 0)   {$D_{T1} O_{T1}$};
    \node[cluster] (3) at (4, 0)   {$O_{T1} D_{T2}$};
    \node[cluster] (4) at (6.25, 0)   {$O_{T1} D_{T2} O_{T2}$};
    \node[cluster] (5) at (9, 0)   {$O_{T1} O_{T2} D_{E1}$};
    \node[cluster] (6) at (12, 0)   {$O_{T2} D_{E1} D_{OS}$};
    \node[cluster] (7) at (12, -2)   {$O_{T2} D_{E1} D_{OS} O_{OS}$};
    \node[cluster] (8) at (7.5, -2)   {$O_{T2} D_{E1} O_{OS} R_{OS}$};
    \node[cluster] (9) at (2.5, -2)   {$O_{T2} D_{E1} O_{OS} R_{OS} D_{DMG}$};
    \node[cluster] (10) at (1, -4)   {$O_{T2} D_{E1} O_{OS} R_{OS} D_{DMG} O_{DMG}$};
    \node[cluster] (11) at (8, -4)   {$O_{T2} D_{E1} O_{OS} R_{OS} O_{DMG} R_{DMG}$};
    \node[cluster] (12) at (12, -6)   {$O_{T2} D_{E1} O_{OS} R_{OS} O_{DMG} R_{DMG} D_{E2}$};
    \node[cluster] (13) at (5, -6)   {$O_{T2} D_{E1} R_{OS} R_{DMG} D_{E2} D_{E3}$};
    \draw[->, thick] (1) -- (2);
    \draw[->, thick] (2) -- (3);
    \draw[->, thick] (3) -- (4);
    \draw[->, thick] (4) -- (5);
    \draw[->, thick] (5) -- (6);
    \draw[->, thick] (6) -- (7);
    \draw[->, thick] (7) -- (8);
    \draw[->, thick] (8) -- (9);
    \draw[->, thick] (9) -- (10);
    \draw[->, thick] (10) -- (11);
    \draw[->, thick] (11) -- (12);
    \draw[->, thick] (12) -- (13);
\end{tikzpicture}
}%
\caption{A rooted junction tree representation of the ID in Fig. \ref{fig:score}}
\label{fig:score_rjt}
\end{figure}

The influence diagram for the problem is presented in Fig. \ref{fig:score}, and a corresponding RJT in Fig. \ref{fig:score_rjt}. For converting the influence diagram into a rooted junction tree, we use the first algorithm from \citet{parmentier2020integer}. The algorithm takes a set of nodes and their information sets, along with a topological order for the nodes and returns a gradual rooted junction tree. A topological order for a graph assigns a unique index to each node, so that for each arc $(i,j) \in A_{ID}$, the index of node $j$ is larger than that of node $i$. 

The first stage of the diagram involves R\&D decisions towards climate sensitivity ($D_{CS}$), damages ($D_{Dmg}$) and technology ($D_{T1}$). If successful, the climate parameter (climate sensitivity and damage exponent) R\&D efforts modify the information structure so that the parametrization is partially revealed in 2050 instead of 2070. This is represented by the nodes $R_{Dmg}$ and $R_{CS}$ and the outcome (success/failure) of the projects by $O_{Dmg}$ and $O_{CS}$. Because the value nodes $v \in N^V$ represent deterministic mappings $s_{I(v)} \to \mathbb{R}$, the value nodes are not explicitly represented in the RJT. Instead, the components of the expected utility can be extracted from the clusters containing $I(v)$, following the ideas used in the computational experiments of \citet{parmentier2020integer}. 

Decisions over emission reductions ($D_{Ei}, i \in \{1,2,3\}$) are made in three stages: in 2030, 2050 and 2070, which represent the medium-term and long-term climate actions. The technological R\&D potentially lowers the costs of emission reductions in 2050 and 2070. We connect our example to the discussion on the feasibility of large-scale deployment of BECCS, which has been a crucial but contested result of many mitigation scenarios \citep{Calvin2021}. The R\&D costs and probabilities for the three levels of BECCS costs are parametrized using expert-elicited estimates in \citet{Baker2015}. These are then reflected in the overall emission reduction costs. It is worth noting the major challenges in long-term technological foresight, which is manifested in the wide spectrum of responses from the experts; but elicitation data is nevertheless useful for illustrating the importance of technological progress. After 2070, the level of climate change is observed based on the chosen emission reductions and the observed branch of climate sensitivity, which then determines the severity of the climate damages along with the observed branch of climate damages.

\subsection{Modifying the influence diagram}

As discussed in Section \ref{section:contributions}, the original formulation is limited to problems with discrete and finite state spaces for all nodes. However, discretizing the emission levels $D_{Ei}$ for $i \in \{1,2,3\}$ would inevitably result in suboptimal solutions, thereby limiting the representability of the naturally continuous nature of the decision variables representing the emission levels. 

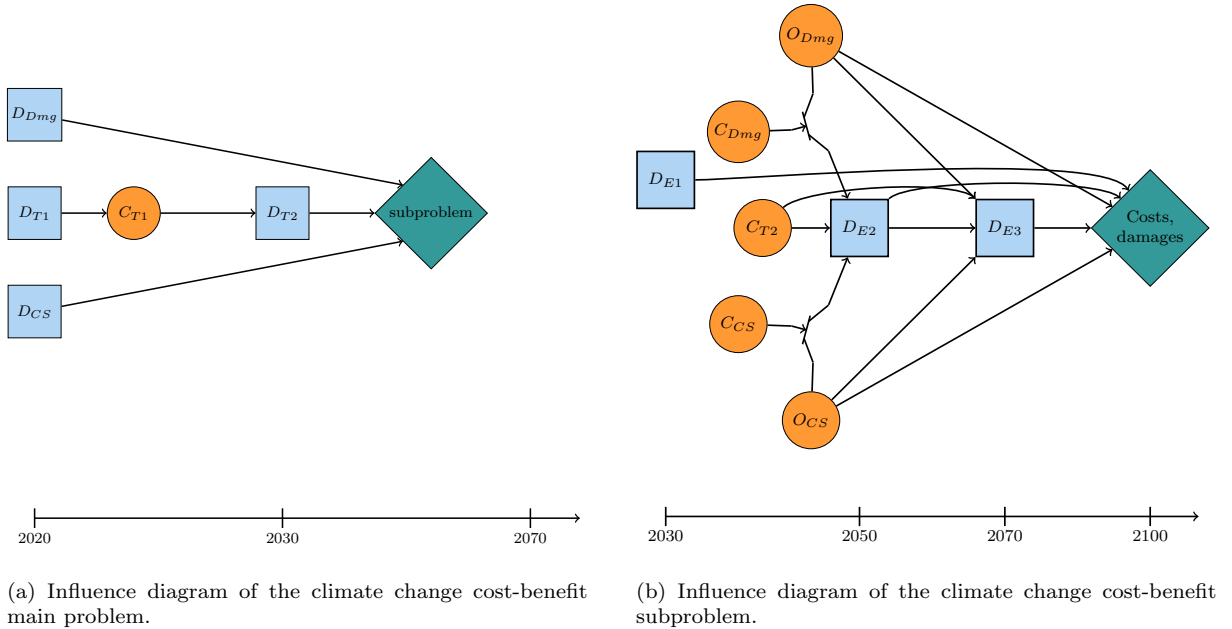
\begin{figure}[!ht]
\centering
     \begin{subfigure}[b]{0.48\textwidth}
        \centering
        \resizebox{\textwidth}{!}{%
        \begin{tikzpicture}
            [decision/.style={fill=blue!40, draw, minimum size=2.5em, inner sep=2pt}, 
            chance/.style={circle, fill=orange!80, draw, minimum size=2.5em, inner sep=2pt},
            value/.style={diamond, fill=teal!80, draw, minimum size=2.5em, inner sep=2pt},
            optimization/.style={ellipse, fill=teal!80, draw, minimum size=2em, inner sep=2pt},
            scale=1.8, font=\scriptsize]
             \node[decision] (D1) at (0, 3)   {$D_{Dmg}$};
             \node[decision] (D2) at (0, 2)   {$D_{T1}$};
             \node[decision] (D3) at (0, 1)   {$D_{CS}$};
             \node[decision] (D5) at (2.5, 2)   {$D_{T2}$};
             \node[chance]   (C2) at (1, 2)   {$O_{T1}$};
             \node[value] (SC) at (4, 2)   {subproblem};
             \draw[->, thick] (D2) -- (C2);
             \draw[->, thick] (C2) -- (D5);
            \draw[->, thick] (D1) -- (SC.135);
            \draw[->, thick] (D3) -- (SC.225);
            \draw[->, thick] (D5) -- (SC);
            \draw[->, thick] (C2) to[out=30,in=150,distance=0.5cm] (SC);
             \draw[->, thick] (0,-1.1) -- (5.5,-1.1);
             \draw[-, thick] (0.0,-1.0) -- (0.0,-1.2);
             \draw[-, thick] (2.5,-1.0) -- (2.5,-1.2);
             \draw[-, thick] (5.0,-1.0) -- (5.0,-1.2);
             \node[] at (0.0,-1.3) {2020};
             \node[] at (2.5,-1.3) {2030};
             \node[] at (5.0,-1.3) {2070};
        \end{tikzpicture}
        }%
        \caption{Influence diagram of the climate change cost-benefit main problem.}
        \label{fig:score_mainproblem}
     \end{subfigure}
     \hfill
     \begin{subfigure}[b]{0.48\textwidth}
        \centering
        \resizebox{\textwidth}{!}{%
        \begin{tikzpicture}
            [decision/.style={fill=blue!40, draw, minimum size=2.5em, inner sep=2pt}, 
            chance/.style={circle, fill=orange!80, draw, minimum size=2.5em, inner sep=2pt},
            value/.style={diamond, fill=teal!80, draw, minimum size=2.5em, inner sep=2pt},
            optimization/.style={ellipse, fill=teal!80, draw, minimum size=2em, inner sep=2pt},
            scale=1.625, font=\scriptsize]
             \node[decision,  thick] (D4) at (2.5, 1.5)   {$D_{E1}$};
             \node[decision,  thick] (D6) at (4.5, 1)   {$D_{E2}$};
             \node[decision,  thick] (D7) at (6.0, 1)   {$D_{E3}$};
             \node[chance]   (C1) at (3.25, 2)   {$O_{Dmg}$};
             \node[chance]   (C3) at (3.25, 0)   {$O_{CS}$};
             \node[chance]   (C4) at (3.5, 1)   {$O_{T2}$};
             \node[chance]   (O1) at (4, 3)   {$R_{Dmg}$};
             \node[chance]   (O2) at (4, -1)   {$R_{CS}$};
             \draw[] (4.1,2.1) pic{dist_arc={int1, , ,rotate=105}};
             \draw[] (4.1,-0.1) pic{dist_arc={int2, , ,rotate=75}};
             \node[value, align=center] (SC) at (7.5, 1)   {Costs, \\damages};
             \draw[->, thick] (O1) -- (D7);
             \draw[->, thick] (O1) -- (SC);
             \draw[->, thick] (O2) -- (D7);
             \draw[->, thick] (O2) -- (SC);
             \draw[->, thick] (D4) to[out=0,in=120,distance=0.5cm] (SC);
             \draw[->, thick] (D6) -- (D7);
             \draw[->, thick] (D6) to[out=45,in=135,distance=0.3cm] (SC);
             \draw[->, thick] (D7) -- (SC);
             
             \draw[->, thick] (C4) -- (D6);
             \draw[->, thick] (C4) to[out=45,in=135,distance=0.3cm] (D7);
             
             \draw[thick] (C1) -- (int1.195);
             \draw[thick] (O1) -- (int1.105);
             \draw[->, thick] (int1.285) -- (D6);
             
             \draw[thick] (C3) -- (int2.165);
             \draw[thick] (O2) -- (int2.255);
             \draw[->, thick] (int2.75) -- (D6);
             
             \draw[->, thick] (2.5,-2) -- (8,-2);
             \draw[-, thick] (2.5,-1.9) -- (2.5,-2.1);
             \draw[-, thick] (4.5,-1.9) -- (4.5,-2.1);
             \draw[-, thick] (6.0,-1.9) -- (6.0,-2.1);
             \draw[-, thick] (7.5,-1.9) -- (7.5,-2.1);
             \node[] at (2.5,-2.2) {2030};
             \node[] at (4.5,-2.2) {2050};
             \node[] at (6.0,-2.2) {2070};
             \node[] at (7.5,-2.2) {2100};
        \end{tikzpicture}
        }%
        \caption{Influence diagram of the climate change cost-benefit subproblem.}
        \label{fig:score_subproblem}
     \end{subfigure}
     \caption{The decomposed climate change cost-benefit problem.}
     \label{fig:score_simplified}
\end{figure}

Following the ideas in Section \ref{section:decomp}, we modify the influence diagram and move the nodes $D_{Ei}$ to a subproblem whose solution becomes the utility value in the main problem. Additionally, we place all chance nodes except $O_{T_1}$ in the subproblem, resulting in the decomposition presented in Fig. \ref{fig:score_simplified}. It should be noted that the choice of which nodes to include in the subproblem is not unique, and for larger models, a significant tradeoff between submodel and main problem sizes might be present, affecting the total solution time. Nevertheless, for this illustrative example, the solution times are small (they are discussed in Section \ref{subsection:score_comp}) making the analysis of secondary importance in our context.

In the resulting subproblem in Fig. \ref{fig:score_subproblem}, there is no Type 1 endogenous uncertainty. Hence,  the subproblem can be modeled as a nonlinear three-stage stochastic programming problem with continuous decision variables. The nonlinearity stems from the SCORE model described in Appendix \ref{sec:app_score}. 

\subsection{Model results}

\begin{figure}[!ht]
    \centering
    \includegraphics[width=\textwidth]{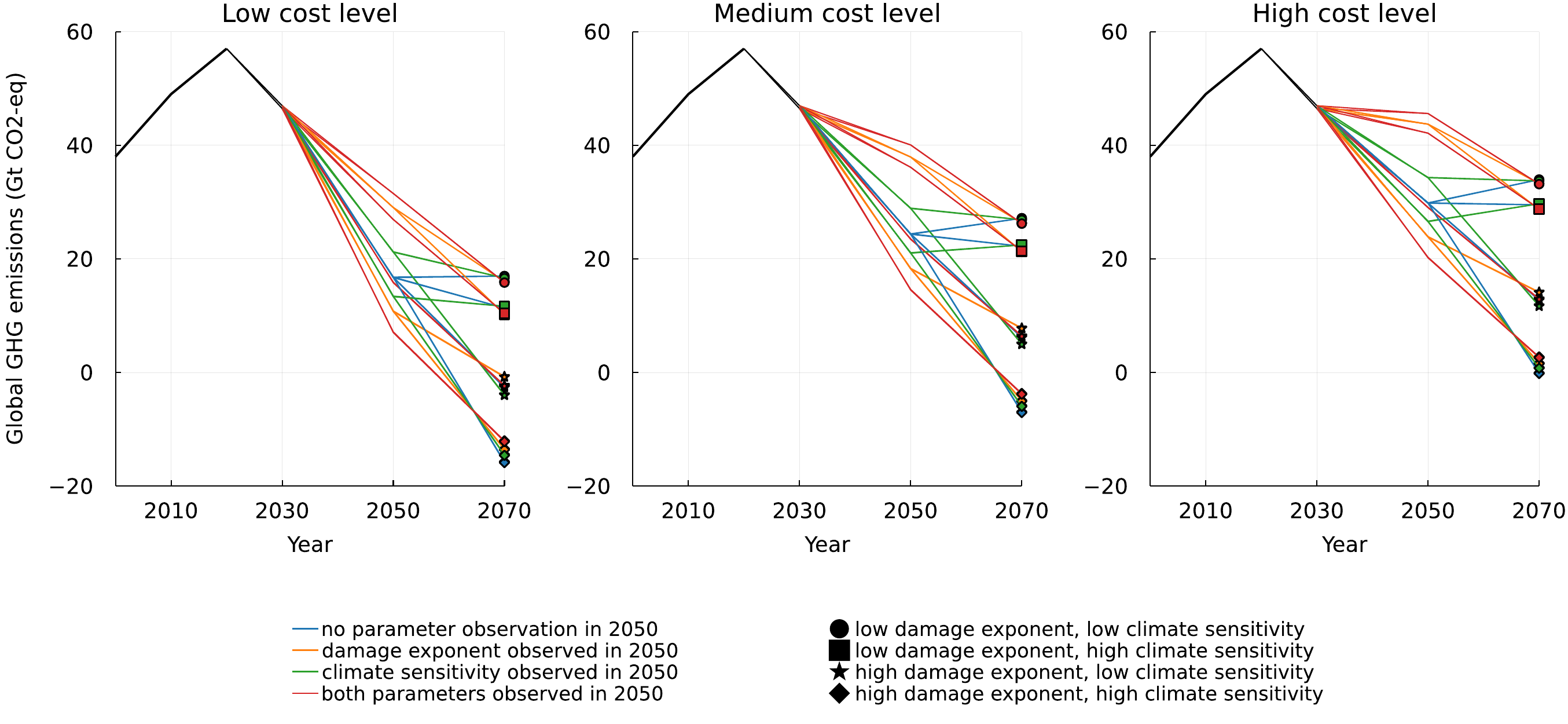}
    \caption{Abatement levels for different outcomes of technology research and parameter observations in 2030.}
    \label{fig:strat_cba}
\end{figure}

Although the primary purpose of this section is to illustrate the computational method, we discuss here some aspects of the model results, as well. The optimal R\&D strategy for this problem is to carry out all R\&D projects. Fig. \ref{fig:strat_cba} presents the emission levels of the optimal mitigation strategy. The effect of the technology R\&D on optimal emission pathways is shown with the three subfigures corresponding to the final R\&D outcome after 2030. Intuitively, successful research and therefore cheaper abatement leads to more abatement. This has also a major impact on the total costs of the optimal strategy: with the low cost curves, the total expected cost is roughly 30\% lower than with medium costs, and with high abatement costs 30\% higher than with the medium cost curve. 

The climatic parameter R\&D, the other endogenous effect in this model has a magnitude smaller effect on the expected costs than the technology R\&D, but the effect on abatement levels is remarkable. The branches after 2030 represent different realizations of the technology research and partial learning of the climate parameters. If the research efforts for both parameters fail (blue lines in Fig. \ref{fig:strat_cba}), the 2050 abatement decisions are made knowing only the outcome of the technology projects, and the four scenarios of partial learning only occur after 2050. Learning the parameters before 2050 results in more dispersed abatement strategies. Finally, the underlying parameter branching is represented with shapes in 2070. It can be seen that the impact of climate sensitivity is considerably smaller than that of the damage exponent. This is in line with the results in \cite{EkholmBaker2021} and shows that the proposed framework could be applied in planning for optimal R\&D pathways.

\subsection{Computational aspects}
\label{subsection:score_comp}

Using the decomposed model presented in Fig. \ref{fig:score_simplified}, the main problem is solved in 0.4 seconds, while solving all of the subproblems takes one second. If we instead solve a problem corresponding to Fig. \ref{fig:score} without decomposition, the framework requires discretizing the abatement levels. Using discrete decision spaces to approximate a continuous variable leads to suboptimal solutions. Furthermore, even with only five levels for each abatement decision, the problem becomes more computationally demanding than the continuous version. The damage costs are calculated in a tenth of a second, but solving the model takes four seconds. As suggested by the experiments in Section \ref{section:experiments}, in larger problems, the discretized problem could quickly become intractable. The decomposed problem instead uses continuous decision variables to represent the abatement decisions, precluding the need for a discretization. As a consequence, it is faster to solve than the rudimentary approximation with five abatement levels.

The discretized model has 2172 constraints and 4332 variables, of which 327 are binary; and the main problem of the decomposed model has 40 constraints and 66 variables, of which 12 are binary. The discretized model is thus two orders of magnitude larger than the decomposed model, but the damage costs in the discretized model are calculated almost instantaneously. However, the trade-off of moving some of the computational burden into the subproblems makes it possible to solve the continuous problem to optimality. Overall, this case study provides an illustrative example of an application where the developments in this paper make it possible to consider settings that are challenging both from modeling and computational standpoints.

\section{Conclusion}
\label{section:conclusions}

In this paper, we propose a framework for modeling and solving Type 3 endogenously uncertain stochastic programming (T3ESP) problems using influence diagrams. The framework is based on the work by \citet{parmentier2020integer}, originally developed for solving decision problems with decision-dependent uncertainties by converting an influence diagram representation of the problem to a mixed-integer programming (MIP) model. To the best of our knowledge, and in line with \citet{hellemo2018}, T3ESP problems have not been previously addressed in the literature.

To make the rooted junction tree (RJT) formulations applicable to T3ESP problems, we show how Type 2 endogenous uncertainties can be modeled by either modifying the underlying influence diagram by adding observation nodes or by incorporating conditional non-anticipativity constraints (C-NACs) to the MIP model. We note that with minimal modifications, the developments presented in this paper can also be applied to extend the Decision Programming framework \citep{Salo2022, hankimaa2023}, as the structure of the two MIP models is sufficiently similar.

In practice, both approaches have their advantages and disadvantages. Adding observation nodes only requires modifying the influence diagram by adding new nodes, while C-NACs are additional constraints that must be added to the MIP model. From the standpoint of a practitioner, it would be simpler to not require explicit modification of the resulting model. On the other hand, if one allows for modifying the MIP model by using C-NACs, this might also allow for representing other parts of the decision process explicitly as added variables or constraints. Further research in this direction is thus relevant.

However, if decision-dependent probability distributions and conditional information revelation are intertwined in the problem structure, for example, by the presence of imperfect conditional observations, one cannot directly employ C-NAC constraints. An example of such a problem would be a version of the climate CBA problem in Section \ref{section:score} where the climate parameter research does not reveal the correct branch, that is, remove one of the extreme parameter values, but instead gives a probability distribution that provides better information than the original. Observation nodes can be used for modeling such uncertainties, while C-NACs require additional modifications to the diagram, further diminishing their benefits. Finally, the MIP models with C-NACs are smaller than the corresponding models using observation nodes, which is in general a desirable feature. In our experiments, however, having smaller MIP models did not translate into improved computational efficiency. Nevertheless, further research is required before more can be concluded in that direction.

Considering alternative modeling paradigms would make the framework suitable for a broader set of problems. Interesting examples of such research ideas are (distributionally) robust optimization and further examination of multi-objective decision-making. For large problems, it might be necessary to improve the computational performance using, e.g., decomposition methods for solving the MIP model \eqref{eq:rjt-obj}-\eqref{eq:rjt-vars}, and solution heuristics. In this paper, we explore the use of influence diagram decomposition \citep[e.g., ][]{lee2021submodel} for improving the computational efficiency of finding maximum expected utility strategies for influence diagrams. Perhaps more interestingly, in the context of our MIP formulations, we also show that such decomposition approaches can allow for solving decision problems with continuous decision variables, significantly improving the general applicability of the formulation.

Our approach is targeted for settings where the T3ESP decision problem can be represented as an influence diagram and associated parameterization. Specifically, our methodological developments enable a rigorous representation of conditionally observed information within influence diagrams by incorporating either observation nodes or conditional non-anticipativity constraints (C-NACs). This, in turn, allows for modeling and solving decision problems with both types of endogenous uncertainty, i.e., decision-dependent probabilities and information structures. Notably, our proposed framework frees the user from being limited to only being able to obtain optimal solutions from their model representation if the diagram complies with particular assumptions. In addition, being based on mathematical programming, the proposed approach benefits from robust off-the-shelf optimization software and precludes the need for time-consuming research, development and implementation of particular solution algorithms. Lastly, while the framework is illustrated with an application in climate change mitigation, the underlying structure applies to a range of decision problems where the timing or availability of information depends on prior decisions, such as technology R\&D portfolio planning, supply chain design, and resource allocation under uncertainty \citep[e.g.,][]{colvin2011r, apap2017, hellemo2018}.

In conclusion, the proposed developments turned the framework sufficiently general to model a challenging example problem with Type 3 endogenous uncertainty. It should also be noted that the influence diagrams are formulated as MILP problems, guaranteeing global optimality of solutions despite the challenging nature of the underlying decision problems.

\section*{Declarations}

\subsection*{Funding}
The work of Herrala and Oliveira was supported by the Research Council of Finland [\textit{Decision Programming: A Stochastic Optimization Framework for Multi-Stage Decision Problems}, grant number 332180]; the work of Ekholm was supported by the Research Council of Finland [\textit{Designing robust climate strategies with Negative Emission Technologies under deep uncertainties and risk accumulation (NETS)}, grant number 331764].

\subsection*{Competing interests}
The authors have no competing interests to declare that are relevant to the content of this article.

\subsection*{Ethics approval and consent to participate}
Not applicable. 

\subsection*{Consent for publication}
All authors mutually agreed that the manuscript should be submitted to Computational Management Science.

\subsection*{Data availability}
Not applicable.

\subsection*{Materials availability}
Not applicable.

\subsection*{Code availability}
The code is available in https://github.com/solliolli/type2-rjt

\subsection*{Author contribution}
\noindent \textbf{Olli Herrala:} Conceptualization, Methodology, Software, Writing - Original Draft, Visualization. \textbf{Tommi Ekholm:} Conceptualization, Investigation, Writing - Review \& Editing. \textbf{Fabricio Oliveira:} Conceptualization, Writing - Review \& Editing, Supervision, Funding acquisition.


\pagebreak
\begin{appendices}
\section{Notation}
\label{sec:app_notation}

\begin{table}[!ht]
\centering
\begin{tabular}{l|l}
symbol & description  \\ \hline
$G_{ID}=(N,A_{ID})$ & an influence diagram consisting of nodes $n \in N$ and arcs $a \in A_{ID}$ \\
$N^C, N^D, N^V$ & sets of chance, decision and value nodes, respectively \\
$I(j)$ & information set of node $j \in N$: $I(j)=\{ i \in N \mid (i,j) \in A_{ID}\}$ \\
$s_j \in S_j$ & state of node $j$ \\
$s_X \in S_X$ & states of nodes in set $X$: $s_X = (s_j)_{j \in X}$ \\  
$G_{RJT}=(V,A_{RJT})$ & a rooted junction tree consisting of clusters $C \in V$ and arcs $(C_i, C_j) \in A_{RJT}$ \\
$Z_j$ & local decision strategy in node $j \in N^D$ mapping $s_{I(j)}$ to a decision $s_j$ \\
$z(s_j \mid s_{I(j)}$ & a binary decision variable with value 1 if and only if $Z_j$ maps $s_{I(j)}$ to $s_j$ \\ 
$\mathbb{P}(X_j=s_j \mid X_{I(j)}=s_{I(j)})$ & conditional probability of node $j \in N^C$ being in state $s_j$ given $s_{I(j)}$ \\
$\mu_{C_j}$ & the joint probability distribution of nodes $n \in C_j$ \\
$u_{C_j}$ & utility function mapping $s_{C_j}$ to a real-valued utility value 
\end{tabular}
\caption{Notation used for influence diagrams and rooted junction trees}
\label{tbl:notation_id_rjt}
\end{table}

\begin{table}[!ht]
\centering
\begin{tabular}{l|l}
symbol & description  \\ \hline
$T_{i,j}$ & distinguishability set, the nodes that are involved \\ & in making node $i \in N^C \cup N^D$ observed in node $j \in N^D$ \\
$F_{i,j}$ & distinguishability condition for observing node $i \in N^C \cup N^D$ in node $j \in N^D$ \\
$a_c \in A_c$ & conditional arcs \\
$I_c(j)$ & conditional information set of node $j$ \\
$f_j^{s_{C_j},s'_{C_j}}$ & a Boolean variable with value 1/True if and only if $s_{C_j}$ and $s'_{C_j}$ are \\ & distinguishable when making decision $j$
\end{tabular}
\caption{Notation used for conditionally observed information}
\label{tbl:notation_type2}
\end{table}
\newpage
\section{N-monitoring parametrization}
\label{sec:app_nmon}
The conditional N-monitoring problem used in the computational experiments is parametrized as follows: 
\begin{itemize}
    \item Node $I_{pre}$ has two states corresponding to low and high initial interest. The probability of the high state is drawn from the uniform distribution $U(0,1)$, and $\mathbb{P}(X_{I_{pre}} = low) = 1-\mathbb{P}(X_{I_{pre}} = high)$.
    \item Nodes $R_i$ have two states each, corresponding to reports of a low and high interest state. The probability of a report in $R_i$ correctly stating a low load is set to $max(x_i,1-x_i)$, where $x_i$ is drawn from $U(0,1)$. Similarly, the probability of correctly stating a high load is $max(y_i,1-y_i)$, with $y_i$ drawn from $U(0,1)$
    \item For the decision nodes $A_i$, the cost $c_i$ of introducing the product in area $i$ is drawn from $U(0,1)$.
    \item The prior probability of the low interest state in $I_{post}$ is $max(x,1-x)$ if $s_{I_{pre}}=low$, and $min(y,1-y)$ if $s_{I_{pre}}=high$. $x$ and $y$ are drawn from $U(0,1)$. This prior probability is then divided by $e^{\sum_{\{i | s_{A_i}="yes"\}}c_i}$, where the exponent is the total cost of introducing the product in each area where it was introduced.
    \item The cost of obtaining a report was set to 10, and the utility $U_{int}$ of the high interest state in $I_{post}$ is 100, while that of the low interest state is 0.
\end{itemize}

\newpage
\section{Cost-benefit model description}
\label{sec:app_score}

The idea in climate change cost-benefit analysis (CBA) is the minimization of emission reduction costs and climate damages. The abatement cost calculation is based on marginal abatement cost curves, as presented in equation \eqref{eq:score_mac}, using numerical estimates from the SCORE model \citep{Ekholm2018}. Coefficients $\alpha$ and $\beta$ are the parameters of the cost curves, $R$ is the total abatement level and $c$ is the marginal cost of abatement. In \eqref{eq:score_abatement_cost}, $C$ is the total cost for abatement level $R$. The subscript $t$ has been omitted for clarity, but parameters $\alpha$ and $\beta$ change between stages due to assumed technological progress.

\begin{align}
    R &= \alpha c^\beta \label{eq:score_mac}\\ 
    \Longrightarrow c &= \left(\frac{R}{\alpha}\right)^{1/\beta}\\
    \Longrightarrow C &= \int_0^R \left(\frac{r}{\alpha}\right)^{1/\beta}dr=\frac{\beta}{1+\beta}\left(\frac{1}{\alpha}\right)^{1/\beta}R^{1+1/\beta}. \label{eq:score_abatement_cost}
\end{align}

The uncertainties considered in this analysis concern 1) the sensitivity of climate to GHG emissions, 2) the severity of climate change damages to society, and 3) the cost of reducing emissions. Decisions can be made to first conduct costly research and development (R\&D) efforts towards each source of uncertainty. For the uncertainties regarding climate sensitivity and damages, a successful R\&D effort results in an earlier revelation of the parametrization, whereas an unsuccessful or no R\&D effort reveals this information later. Similar models of R\&D pipeline optimization under (Type 2) endogenous uncertainty are considered in \citet{colvin2011r}.

Departing from the predetermined technological progress that was assumed by \cite{Ekholm2018}, the parameter $\alpha$ in the model depends here on the result of technological R\&D, as presented in Fig. \ref{fig:score}. We consider three levels of R\&D effort in 2020 and 2030, which can then lead to three possible levels of MAC curves for years 2050 and 2070. For a related discussion, we refer the reader to \cite{rathi2022capacity}, who consider endogenous technology learning on power generation.

For the effect of R\&D efforts on bioenergy and carbon capture and storage (CCS) costs, we used expert elicited estimates from \citet{Baker2015}. To convert these into emission reduction costs, we assume a coal power plant as a baseline and calculate the additional costs from bioenergy with carbon capture and storage (BECCS) relative to the amount of reduced emissions by switching from coal to BECCS. Both plants were assumed to have a lifetime of 30 years and operate at 80\% capacity on average. The coal power plant was assumed to have a 40\% efficiency and produce 885 tonnes of CO$_2$ per GWh of electricity. The generation cost was assumed to be 50\$/MWh. Costs were discounted at 5\% rate. 

Compared to coal, BECCS accrues additional costs per generated unit of electricity from the higher cost of biofuel and lower efficiency, and the additional investment to CCS and loss of efficiency from using some of the generated electricity in the carbon capture process. 
\citet{Baker2015} presented probability distributions for these parameters following three different levels of R\&D efforts (low, medium and high). To calculate the cost differential to coal power plant, we performed a Monte Carlo sampling of these four parameters, separately for each R\&D level, which was then compared to the amount of reduced emissions per generated unit of electricity. This yields a distribution of emission reduction costs for BECCS for each R\&D level.

We generalize the impact of R\&D on BECCS's emission reduction cost to the overall marginal abatement cost (MAC) curve. This is obviously a simplification, but nevertheless reflects the major role that BECCS might have in decarbonizing the economy \citep[e.g. ][]{Fuss2018,Rogelj2018}. We take the high MAC from \citet{Ekholm2018} as the starting point and define two MAC curves that are proportionally scaled down from the high MAC.

The low R\&D level yields an average emission reduction cost of around 100 \$/t. We set three bins for three cost levels: high costs correspond to above 75 \$/t, medium costs are between 25 and 75 \$/t, and low costs are below 25 \$/t. The probabilities of achieving high, medium or low emission reduction costs are then estimated from the Monte Carlo sampling for each R\&D level, and presented in Table \ref{tbl:MAC_prob}. With medium and high R\&D effort, the average cost in the low cost bin is around 20 \$/t. Therefore we assign the reduction in the MAC as 50\% for medium costs and 80\% for low costs. The corresponding parameters are listed in Table \ref{tbl:score_ab}. These are still within the range of costs used in \citet{Ekholm2018}, where the low-cost MAC yielded the same emission reductions than the high-cost MAC with approximately 90\% lower costs in 2050.

\begin{table}
\centering
\caption{Probabilities for different abatement costs (rows) with different levels of R\&D effort (columns).}
\label{tbl:MAC_prob}
\begin{tabular}{lccc}
            & Low R\&D  & Medium R\&D   & High R\&D   \\ \hline
High costs  & 73 \%     & 31 \%         & 9 \%  \\
Medium costs& 27 \%     & 64 \%         & 73 \% \\
Low costs   & 0 \%      &  5 \%         & 18 \%
\end{tabular}
\end{table}

\begin{table}
\centering
\caption{Cost curve parametrization}
\label{tbl:score_ab}
\begin{tabular}{lcccc}
year & $\alpha_{high}$  & $\alpha_{medium}$ & $\alpha_{low}$    & $\beta$   \\ \hline
2030 & 3.57             & 3.57              & 3.57              & 0.340     \\
2050 & 11.2             & 13.3              & 16.7              & 0.250     \\
2070 & 21.1             & 24.3              & 29.3              & 0.203
\end{tabular}
\end{table}

The climate damage cost calculation is from DICE \citep{Nordhaus2017}. The damage function is presented in \eqref{eq:score_dice}, where $Y(t)$ is the world gross economic output at time $t$, $a$ is a scaling parameter and $b$ is the damage exponent. While climate change and the abatement decisions have an effect on the economic output, the effect is assumed small and $Y(t)$ is defined exogenously in SCORE.

\begin{equation}
    D(t,\Delta T) = Y(t)a\Delta T^b. \label{eq:score_dice}
\end{equation}

Finally, the temperature change $\Delta T$ is approximated with \eqref{eq:score_deltat}, where $c$ is the climate sensitivity (the temperature increase from doubling of CO$_2$ emissions), $M$ is the sum of emissions in 2030-2070 and $k_i$ are coefficients. 

\begin{equation}
    \Delta T = k_1 c M + k_2 c + k_3 M + k_4. \label{eq:score_deltat}
\end{equation}

In SCORE, both the DICE damage parameter and the climate sensitivity are uncertain with three options, low, medium and high, as presented in Table \ref{tbl:param_values}. The uncertainty is revealed in two steps in a binomial lattice, which is a simplified version of the lattice employed in \citet{Ekholm2018}. First, between 2050 and 2070, one of the extreme alternatives is removed from both uncertainties, that is, for both parameters, we know either that the value is not high or that it is not low. Then, after 2070, we learn the actual value. 

The implementation here combines influence diagrams and multi-stage stochastic programming (MSSP) in a way that the underlying branching probabilities are used as the probabilities of the observations $R_{Dmg}$ and $R_{CS}$ in Fig. \ref{fig:score}. For the damage exponent, all branching probabilities are 50\%. The observation $R_{Dmg}$ thus has a 50\% probability of removing either the high or low value. Depending on the branch, the low or high value then has a 50\% probability in the later branching, with the other 50\% for the medium value. This makes the medium branch have a 50\% probability in total, while the two extreme values both have a 25\% probability. For the climate sensitivity, the first branching is with a 50\% probability for both branches. However, the second branching is different. If the high sensitivity is excluded in the first branch, there is a 21\% conditional probability of the low branch in the second branching, meaning a 10.5\% total probability for the low sensitivity. Similarly, there is a 23\% conditional probability of high damages in the other branch, resulting in a 11.5\% probability for the high sensitivity. The remaining 78\% is the final probability of medium sensitivity. 

\begin{table}
\centering
\caption{Climate sensitivity and damage exponent values}
\label{tbl:param_values}
\begin{tabular}{lccc}
        & Climate sensitivity   & Damage exponent   \\ \hline
High    & 6                     & 4                 \\
Medium  & 3                     & 2                 \\
Low     & 1.5                   & 1
\end{tabular}
\end{table}

This uncertain process is modeled by means of a multi-stage stochastic programming problem, where new information is obtained gradually. It is possible to perform research on these parameters. If the research succeeds, one of the extreme values is excluded already before 2050, revealing the first branching in the observation process earlier than without or with failed research. The observation of the actual parameter value (the second branching) still happens after 2070, after all abatement decisions have been made. The total cost we aim to minimize is then a discounted sum of research costs, abatement costs \eqref{eq:score_abatement_cost} for years 2030, 2050 and 2070, and damage costs \eqref{eq:score_dice}.
\newpage
\section{MIP formulation of the 2-monitoring problem}
\label{app:nmon}

\begin{align}
    \max &\sum_{s_{I_{pre}} \in S_{I_{pre}}, s_{A_1} \in S_{A_1}, s_{V_1} \in S_{V_1}} \mu_{C_{V_1}}(s_{I_{pre}},s_{A_1},s_{V_1}) u_{V_1}(s_{I_{pre}},s_{A_1},s_{V_1}) \nonumber \\ 
    & \quad + \sum_{s_{I_{pre}} \in S_{I_{pre}}, s_{A_2} \in S_{A_2}, s_{V_2} \in S_{V_2}} \mu_{C_{V_2}}(s_{I_{pre}},s_{A_2},s_{V_2}) u_{V_2}(s_{I_{pre}},s_{A_2},s_{V_2}) \nonumber \\ 
    & \quad + \sum_{s_{I_{post}} \in S_{I_{post}}, s_{U_{int}} \in S_{U_{int}}} \mu_{C_{U_{int}}}(s_{I_{post}},s_{U_{int}}) u_{U_{int}}(s_{I_{post}},s_{U_{int}}) \\
    \text{s.t. } & \sum_{s_{I_{pre}} \in S_{I_{pre}}} \mu_{C_{I_{pre}}}(s_{I_{pre}}) = 1, \\
    & \mu_{C_{I_{pre}}}(s_{I_{pre}}) = \mathbb{P}(X_{I_{pre}}=s_{I_{pre}}), \ \forall s_{{I_{pre}}} \in S_{{I_{pre}}} \\
    & \mu_{C_{I_{pre}}}(s_{{I_{pre}}}) \ge 0, \ \forall s_{{I_{pre}}} \in S_{{I_{pre}}}\\
    & \sum_{s_{I_{pre}} \in S_{I_{pre}}, s_{R_1} \in S_{R_1}} \mu_{C_{R_1}}(s_{I_{pre}},s_{R_1}) = 1, \\
    & \mu_{C_{I_{pre}}}(s^*_{{I_{pre}}}) = \sum_{s_{R_1} \in S_{R_1}} \mu_{C_{R_1}}(s^*_{I_{pre}}, s_{R_1}), \ \forall  s^*_{{I_{pre}}} \in S_{{I_{pre}}} \\
    & \mu_{C_{R_1}}(s_{I_{pre}},s_{R_1}) = \mu_{\overline{C}_{R_1}}(s_{I_{pre}}) \mathbb{P}(X_{R_1}=s_{R_1} \mid X_{{I_{pre}}}=s_{{I_{pre}}}), \ \forall s_{{I_{pre}}} \in S_{{I_{pre}}}, s_{R_1} \in S_{R_1} \\
    & \mu_{C_{R_1}}(s_{I_{pre}},s_{R_1}) \ge 0, \ \forall s_{I_{pre}} \in S_{I_{pre}}, s_{R_1} \in S_{R_1}\\
    & \sum_{s_{I_{pre}} \in S_{I_{pre}}, s_{R_1} \in S_{R_1}, s_{A_1} \in S_{A_1}} \mu_{C_{A_1}}(s_{I_{pre}},s_{R_1},s_{A_1}) = 1, \\
    & \mu_{C_{R_1}}(s^*_{{I_{pre}}}, s^*_{R_1}) = \sum_{s_{A_1} \in S_{A_1}} \mu_{C_{A_1}}(s^*_{I_{pre}}, s^*_{R_1}, s_{A_1}), \ \forall s^*_{{I_{pre}}} \in S_{{I_{pre}}}, s^*_{R_1} \in S^*_{R_1}\\
    & \mu_{C_{A_1}}(s_{I_{pre}}, s_{R_1}, s_{A_1}) = \mu_{\overline{C}_{A_1}}(s_{I_{pre}}, s_{R_1})z(s_{A_1} \mid s_{R_1}), \ \forall s_{I_{pre}} \in S_{I_{pre}}, s_{R_1} \in S_{R_1}, s_{A_1} \in S_{A_1}\\
    & \mu_{C_{A_1}}(s_{I_{pre}},s_{R_1},s_{A_1}) \ge 0, \ \forall s_{I_{pre}} \in S_{I_{pre}}, s_{R_1} \in S_{R_1}, s_{A_1} \in S_{A_1}\\
    & z(s_{A_1} \mid s_{R_1}) \in \{0,1\}, \ \forall s_{R_1} \in S_{R_1}, s_{A_1} \in S_{A_1} \\
    & \sum_{s_{A_1} \in S_{A_1}, s_{V_1} \in S_{V_1}} \mu_{C_{V_1}}(s_{A_1},s_{V_1}) = 1, \\
    & \sum_{s_{{I_{pre}}} \in S_{{I_{pre}}}, s_{R_1} \in S_{R_1}} \mu_{C_{A_1}}(s_{{I_{pre}}}, s_{R_1}, s^*_{A_1}) = \sum_{s_{V_1} \in S_{V_1}} \mu_{C_{V_1}}(s^*_{A_1}, s_{V_1}), \ \forall s^*_{A_1} \in S_{A_1}\\
    & \mu_{C_{V_1}}(s_{A_1},s_{V_1}) = \mu_{\overline{C}_{V_1}}(s_{A_1}) \mathbb{P}(X_{V_1}=s_{V_1} \mid X_{A_1}=s_{A_1}), \ \forall s_{A_1} \in S_{A_1}, s_{V_1} \in S_{V_1} \\
    & \mu_{C_{V_1}}(s_{A_1},s_{V_1}) \ge 0, \ \forall s_{A_1} \in S_{A_1}, s_{V_1} \in S_{V_1}\\
    & \quad ... \nonumber  
\end{align}

\end{appendices}

\bibliography{library}

\end{document}